\documentstyle[12pt]{article}                  
\input{amssymb.sty}
\makeatletter
\def\nothing#1{}
\newdimen\earraycolsep
\renewcommand{\thetable}{\arabic{table}}
\renewcommand{\thefigure}{\arabic{figure}}
\renewcommand{\title}[1]{%
  \vspace*{120\p@}%
  {\parindent \z@ \raggedright \reset@font
    \bfseries #1\par
    \nobreak
    \vskip 36\p@
  }}
\def\author#1{{\pretolerance=10000 \raggedright \advance \leftskip by 1in
\noindent #1 \vskip 1pc}}
\def\affiliation#1{{\advance\leftskip by 1in \noindent #1 \vskip -1pc}}
\def\refnote#1{{$^{\hbox{\scriptsize #1}}$}}
\def\affnote#1{\llap{$^{\hbox{\scriptsize #1}}$}}

\renewcommand\section{\@startsection{section}{1}{\z@}{2pc \@plus
      1ex minus .2ex}{1pc \@plus .2ex}{\reset@font
      \normalsize\bfseries\noindent
      {\addtocounter{section}{1}}\arabic{section}\
      {\setcounter{subsection}{0}
      \setcounter{subsubsection}{0}\setcounter{equation}{0}} }}
\renewcommand\subsection{\@startsection{subsection}{2}{\z@}{1pc \@plus 1ex
    minus.2ex}{1pc \@plus .2ex}
    {\reset@font\normalsize\bfseries
    \noindent{\addtocounter{subsection}{1}}%
    {\setcounter{subsubsection}{0}}\arabic{section}.\arabic{subsection}\ }}
\renewcommand\subsubsection{\@startsection{subsubsection}{3}{\parindent}
        {1pc \@plus 1ex minus.2ex}{-0.5em}{\reset@font\normalsize\bfseries%
        {\addtocounter{subsubsection}{1}} \hspace*{.6cm}
        \arabic{section}.\arabic{subsection}.\arabic{subsubsection}
        \hspace*{-7mm}}}

\def\AmS{{\protect\the\textfont2%
        A\kern-.1667em\lower.5ex\hbox{M}\kern-.125emS}}

\def\p@LaTeX{{\family{times}\series{m}\shape{n}\selectfont
L\kern-.36em\raise.3ex\hbox{\scriptsize A}\kern-.15em
T\kern-.1667em\lower.7ex\hbox{E}\kern-.125emX}}

\newlength{\colwidth}

\def\@oddhead{\hfil}
\def\@evenhead{\hfil}
\def\@oddfoot{{\bfseries\hfil\thepage}}
\def\@evenfoot{{\bfseries\thepage\hfil}}
\def\fnum@figure{\footnotesize\raggedright{\bfseries \figurename~\thefigure.}}
\def\fnum@table{\normalsize\raggedright{\bfseries \tablename~\thetable.}}
\long\def\@makecaption#1#2{\vskip 10\p@ {#1 #2\par}}
\long\def\@makefntext#1{\setbox0=\hbox{$\m@th^{\@thefnmark}$}\noindent
\hangindent=\wd0 \box0 #1}
\newbox\@atbox
\long\def\atable#1#2#3{\begin{table}[tbp]\centering\footnotesize
\setbox\@atbox\hbox{#2}
\parbox{\wd\@atbox}{\caption{#1}}\par\smallskip
#2
\par\smallskip\parbox{\wd\@atbox}{\raggedright #3}
\end{table}}

\newtheorem{theorem}{Theorem}

\newtheorem{proposition}[theorem]{Proposition}

\newtheorem{remark}[theorem]{Remark}

\newcommand{\ra}{\rightarrow}
\newcommand{\fl}{\forall}
\newcommand{\wt}{\widetilde}

\newcommand{\s}{\sigma}
\newcommand{\D}{\Delta}

\newcommand{\Zb}{\mathbb{Z}}
\newcommand{\Uc}{\mathcal{U}}
\newcommand{\ot}{\otimes}

\newcommand{\Hc}{\mathcal{H}}
\newcommand{\g}{\gamma}
\newcommand{\vp}{\varphi}
\newcommand{\ve}{\varepsilon}
\newcommand{\Cb}{\mathbb{C}}

\newcommand{\FG}{{\mathfrak{g}}}

\def\Cb{{\mathbb C}}

\def\Rb{{\mathbb R}}

\def\Zb{{\mathbb Z}}

\def\Ac{{\cal A}}

\def\Ec{{\cal E}}
\def\Fc{{\cal F}}

\def\Hc{{\cal H}}

\def\Lc{{\cal L}}
\def\Nc{{\cal N}}
\def\Rc{{\cal R}}
\def\Sc{{\cal S}}
\def\Uc{{\cal U}}
\def\Vc{{\cal V}}

\def\a{\alpha}

\def\d{\delta}
\def\lb{\lambda}
\def\g{\gamma}
\def\om{\omega}
\def\s{\sigma}

\def\ve{\varepsilon}
\def\vp{\varphi}

\def\D{\Delta}
\def\G{\Gamma}

\def\fl{\forall}
\def\ify{\infty}
\def\lgl{\langle}
\def\nb{\nabla}
\def\op{\oplus}
\def\ot{\otimes}
\def\part{\partial}
\def\rgl{\rangle}
\def\sbs{\subset}
\def\semi{>\!\!\!\lhd}

\def\ra{\rightarrow}

\def\text{\hbox}

\def\Ker{\mathop{\rm Ker}\nolimits}

\def\Sign{\mathop{\rm Sign}\nolimits}

\def\Trace{\mathop{\rm Trace}\nolimits}

\def\boxit#1#2{\setbox1=\hbox{\kern#1{#2}\kern#1}%
\dimen1=\ht1 \advance\dimen1 by #1 \dimen2=\dp1 \advance\dimen2 by #1
\setbox1=\hbox{\vrule height\dimen1 depth\dimen2\box1\vrule}%
\setbox1=\vbox{\hrule\box1\hrule}%
\advance\dimen1 by .4pt \ht1=\dimen1
\advance\dimen2 by .4pt \dp1=\dimen2 \box1\relax}

\def\@nbibitem#1{\noindent \hangindent=2pc \hangafter=1
\refstepcounter{enumi}\hbox to 2pc{\arabic{enumi}.\hfil}%
\immediate\write\@auxout{\string\bibcite{#1}{\arabic{enumi}}}}
\def\numbibliography{%
\section*{REFERENCES}%
\bgroup\footnotesize
\setcounter{enumi}{0}%
\def\newblock{\hskip .11em plus.33em minus.07em}%
\let\bibitem\@nbibitem}
\def\endnumbibliography{\par\egroup}
\makeatother
\begin{document}

\begin{center}
{ \bf CYCLIC COHOMOLOGY AND HOPF SYMMETRY}
\end{center}
\medskip
\author{\bf Alain CONNES\refnote{1} and 
    Henri MOSCOVICI\refnote{2}\footnote{The second author has 
    conducted this work during employment by the Clay Mathematics 
    Institute, as a CMI Scholar in residence at Harvard University. 
    This material is based in part 
    upon research supported by the National Science Foundation
    under award no. DMS-9706886.}}

\affiliation{\affnote{1}  Coll\`ege de France,
3, rue Ulm,
75005 Paris\\
and\\
I.H.E.S., 35, route de Chartres, 91440 Bures-sur-Yvette \\
\affnote{2} Department of Mathematics,
The Ohio State University \\
231 W. 18th Avenue, Columbus, OH 43210, USA  
}
\medskip

\begin{abstract}
Cyclic cohomology has been recently adapted to the 
treatment of Hopf symmetry in noncommutative geometry.
The resulting theory of characteristic classes 
for Hopf algebras and their actions on algebras allows to 
expand the range of applications of cyclic cohomology. It is
the goal of the present paper to illustrate these recent
developments, with special emphasis on the application to
transverse index theory,
and point towards future directions. In particular,
we highlight the remarkable accord between
our framework for cyclic cohomology of Hopf algebras on one hand
and both the algebraic as well as the analytic 
theory of quantum groups on the other, manifest in the
construction of the modular square.
\end{abstract}
\vspace{1cm}

\centerline{\bf Introduction}
\medskip

Cyclic cohomology of noncommutative algebras is playing 
in noncommutative geometry a similar r\^ole to that of de Rham cohomology 
in differential topology \cite{Co}.  
In \cite{CM2} and \cite{CM3}, cyclic cohomology has been adapted
to Hopf algebras and their
actions on algebras, which are analogous to the Lie group/algebra
actions on manifolds and embody a natural notion of symmetry
in noncommutative geometry. The resulting theory of
characteristic classes for Hopf actions allows in turn to widen
the scope of applications of cyclic cohomology to index theory.
It is the goal of the present paper to review these recent 
developments and point towards future directions.
\smallskip

The contents of the paper are as follows. In \S 1 we 
recall the basic notation pertaining to the cyclic theory. The 
adaptation of cyclic cohomology to Hopf algebras and Hopf actions
is reviewed in \S 2, where we also discuss the relationship with
Lie group/algebra cohomology. \S 3 deals with the geometric Hopf algebras
arising in transverse differential geometry and their application
to transverse index theory. Finally, \S 4 illustrates the remarkable
agreement between our framework for cyclic cohomology of Hopf algebras
and both the algebraic as well as the analytic theory of quantum groups.
\vspace{1cm}

\section{Cyclic cohomology}
 
Cyclic cohomology has first appeared as a cohomology theory for 
algebras (\cite{C0}, \cite{Co0}, \cite{Ts}). 
In its simplest form, the cyclic cohomology $ HC^{*} (\Ac)$ 
of an algebra  $\Ac$ (over $\Rb$ or $\Cb$ in what 
follows) is the cohomology of 
the cochain complex $\{ C_{\lambda}^{*} (\Ac), \, b \} $, where
$ C_{\lambda}^{n} (\Ac)$, $n \geq 0$, consists 
of the ($n+1$)-linear forms $\vp$ on $\Ac$ satisfying the cyclicity condition
\begin{equation} 
\varphi (a^0 ,a^1 ,...,a^n) = (-1)^n \varphi (a^1 , a^2 ,...,a^0) , \quad 
a^0 ,a^1 ,...,a^n \in \Ac \label{eq:cyc}
\end{equation} 
and the coboundary operator is given by 
\begin{eqnarray}
( b\varphi) (a^0 ,\ldots ,a^{n+1}) = 
\sum_{j=0}^n \, (-1)^j \, \varphi (a^0 ,\ldots ,a^j a^{j+1}, \ldots
,a^{n+1})  \nonumber \\
  + (-1)^{n+1} \, \varphi (a^{n+1} a^0 ,a^1 ,\ldots ,a^{n}) .
 \label{eq:cob}
\end{eqnarray}
When the algebra $\Ac$ comes equipped with a locally convex topology
for which the product is continuous, the above complex is replaced
by its topological version: 
$ C_{\lambda}^{n} (\Ac)$ then consists of all continuous
($n+1$)-linear form on $\Ac$ satisfying (\ref{eq:cyc}).

Cyclic cohomology provides numerical 
invariants of K-theory classes as follows (\cite{C2}).
Given an n-dimensional cyclic cocycle $\varphi$
on $\Ac$, $n$ even, the scalar 
\begin{equation}
\varphi \otimes {\rm Tr} \, (E,E,...,E) \label{eq:pairing}
\end{equation}
is invariant 
under homotopy for idempotents 
$$E^{2} = E \in M_N (\Ac) = 
\Ac \otimes  M_N (\Cb) .$$
In the above formula, $\varphi  \otimes {\rm Tr}$ is the extension of $\varphi$
to $M_N (\Ac)$, using the standard trace 
${\rm Tr}$ on $ M_N (\Cb) $:
$$\varphi \otimes {\rm Tr} \,
(a^0 \otimes \mu^0 , a^1 \otimes \mu^1 , \ldots , a^n \otimes \mu^n) 
= \varphi (a^0 , a^1 , \ldots , a^n) \, {\rm Tr} (\mu^0 \mu^1 \cdots \mu^n). 
$$
This defines a pairing $ \langle [\varphi] , [E] \rangle$ between cyclic 
cohomology and $K$--theory, which extends to the general noncommutative
framework the Chern-Weil construction of
characteristic classes of vector bundles.

Indeed,
if $\Ac = C^{\infty} (M)$ for a closed manifold $M$ and 
$$
\varphi (f^0 ,f^1 ,...,f^n) = \langle \Phi , \, f^0 df^1 \wedge df^2\wedge 
\ldots \wedge df^n 
\rangle , \quad f^0 ,f^1 ,...,f^n \in \Ac ,
$$
where $\Phi $ is an n-dimensional closed de Rham current on $M$, then up to 
normalization the invariant defined by (\ref{eq:pairing}) is equal to 
$$
\langle \Phi , \, ch^{*}€ (\Ec) \rangle ;
$$
here $ch^{*}€ (\Ec)$ denotes the Chern character of the rank $N$
vector bundle $\Ec$  
on $M$ whose fiber at $x \in M$ is the range of $E(x) \in M_N (\Cb)$. 

Note that, in the above example,
$$
\varphi (f^{\sigma (0)} , f^{\sigma (1)} ,..., f^{\sigma (n)}) = 
\varepsilon(\sigma) \, \varphi (f^0 , f^1 ,..., f^n) \, ,
$$
for any permutation $\sigma$ of the set $[n] = \{ 0,1,...,n \}$,
with signature $\varepsilon(\sigma) $. However, the extension  
$\varphi \otimes {\rm Tr}$ to
$M_N (\Ac)$, used in the pairing formula (\ref{eq:pairing}), 
retains only the {\it cyclic} invariance.

A simple but very useful class of
examples of cyclic cocycles on a noncommutative algebra is obtained 
from group cohomology (\cite{C4}, \cite{Bu}, \cite{CM0}), as
follows. Let $\Gamma$ be an arbitrary group and let  $\Ac = \Cb \Gamma$
be its group ring. Then
any {\it normalized group cocycle}  $c \in Z^n (\Gamma ,\Cb)$ ,
representing an arbitrary cohomology class 
$[c] \in H^* (B\Gamma) = H^* (\Gamma)$, 
gives rise to a cyclic
cocycle $\varphi_c$ on the algebra $\Ac$ by means of the formula
\begin{equation}
\varphi_c (g_0 , g_1 , \ldots , g_n) = \left\{ \matrix{ 0 \ \hbox{if} \ g_0
\ldots g_n \not= 1 \cr c(g_1 , \ldots , g_n) \ \hbox{if} \ g_0 \ldots g_n =1
\cr} \right. \nonumber
\end{equation} 
extended by linearity to $\Cb \Gamma$.
\smallskip

In a dual fashion, one defines the cyclic homology $ HC_{*} (\Ac)$
of an algebra $\Ac$
as the homology of the chain complex $\{ C^{\lambda}_{*} (\Ac), \, b \} $
consisting of the
coinvariants under cyclic permutations of
the tensor powers 
of $\Ac$, and with the boundary operator $b$ obtained by transposing
the coboundary formula (\ref{eq:cob}). Then the pairing between cyclic
cohomology and $K$--theory 
(\ref{eq:pairing}) factors through the natural pairing between
cohomology and homology, i.e.
\begin{equation}
\varphi \otimes {\rm Tr} \, (E,E,...,E) = \, \langle \varphi, \, ch (E) 
\rangle \, ,
\label{eq:hpairing}
\end{equation}
where, again up to normalization,
\begin{equation} \label{chn}
ch_{n}€ (E) = \, E \ot  ...  \ot E \, \quad (n+1 \quad {\rm times}) 
\label{eq:hch}
\end{equation}
represents the {\it Chern character} in $ HC_{n} (\Ac)$, for $n$ 
even,
of the $K$--theory class $[E] \in
K_{0}€(\Ac)$.
\medskip 

The cyclic cohomology of an (unital) algebra $\Ac$ 
has an equivalent description,
in terms of the bicomplex 
$( C C^{*, *} (\Ac), \, b , \, B ) ,$
defined as follows. With $C^{n} (\Ac)$ denoting the linear space
of ($n+1$)-linear forms on $A$, set
\begin{equation}\label{abB}
\matrix{ 
&C C^{p, q} (\Ac) = C^{q-p} (\Ac), \quad q \geq p , \cr \cr
& C C^{p, q} (\Ac) = 0 ,  \quad q < p \, . \hfill \cr
}
\end{equation}
The vertical operator 
$b: C^{n} (\Ac) \ra C^{n+1} (\Ac) $
is defined as
\begin{eqnarray}
\lefteqn{(b\vp) (a^0 ,\ldots ,a^{n+1}) =} \nonumber \\
& \sum_0^n (-1)^j \, \vp (a^0 ,\ldots ,a^j a^{j+1}, \ldots
,a^{n+1}) \nonumber \\
& + (-1)^{n+1} \, \vp (a^{n+1} a^0 ,a^1 ,\ldots ,
a^{n})  .                                               \label{eq:ab}
\end{eqnarray}
The horizontal operator
$B: C^{n} (\Ac) \ra C^{n-1} (\Ac)$ 
is defined by the formula 
$$ B = N B_0 , $$
where
\begin{eqnarray}
\lefteqn{B_0 \, \vp (a^0 ,\ldots ,a^{n-1}) = \vp
(1,a^0 ,\ldots ,a^{n-1}) - (-1)^n \, \vp (a^0 ,\ldots
,a^{n-1} ,1)} \, , \nonumber \\
& \qquad (N\psi) (a^0 ,\ldots ,a^{n}) = \sum_0^{n}
(-1)^{n j} \, \psi (a^j ,a^{j+1}, \ldots ,a^{j-1}) \,
.                                                \label{eq:aB}
\end{eqnarray}
Then $H C ^{*}€(\Ac)$ is the cohomology of the {\it first
quadrant} total complex
$ ( TC^{*}€(\Ac),\, b+B )$, formed as follows:
\begin{equation}
    TC^{n}€(\Ac) = \sum_{p=0}^{n}€€C C^{p, n-p} (\Ac) . \label{T}
\end{equation}
On the other hand, the cohomology of the {\it full direct sum} 
total complex $( T C^{\Sigma}€_{*}€(\Ac),\,  b+B )$, formed
by taking direct sums as follows:
\begin{equation}
    T C^{\Sigma}€_{n}€(\Ac) = \sum_{s}€€ \, C C^{p, n-p} (\Ac) ,
    \label{ST}
\end{equation}
gives the ( $\Zb /2$--graded) {\it periodic cyclic cohomology} groups
$H C^{*}€_{\rm per}€(\Ac)$.

There is a dual description for the cyclic homology of $\Ac$, in
terms of the dual bicomplex
$( C C_{*, *} (\Ac), \, b , \, B )$, with
$ C_{n}€(\Ac) = \Ac^{\ot n+1}€$ 
and the boundary operators $b$, $B$ obtained by transposing
the corresponding coboundaries. The 
periodic cyclic homology groups
$H C_{*}€^{\rm per}€(\Ac)$ are obtained from the {\it full 
product} total complex $( T C_{\Pi}€^{*}€(\Ac),\,  b+B )$,
formed by taking direct products as follows:
\begin{equation}
   T C_{\Pi}€^{n}€(\Ac) = \Pi_{p} \,  C C_{p, n-p} (\Ac) .
    \label{PT}
\end{equation}
The Chern character of an idempotent $e^{2}€=e \in \Ac$ is given
in this picture by the periodic cycle 
$ ( ch_{n}€ (e) )_{n=2,4,\ldots}€ \,$, with components:
\begin{equation} \label{cch}
 ch_{0}€(e) = e, \quad 
 ch_{2k}€(e) = (-1)^{k}€{(2k)! \over k!}  (e^{\ot 2k+1}€ - {1 \over 2}
    \ot e^{\ot 2k}€ ) , \quad k \geq 1 .
\end{equation}   
\smallskip

The functors $ HC^{0}€$ and $ HC_{0}$ from the category of algebras
to the category of vector spaces have clear intrinsic meaning: the first
assigns to an algebra $\Ac$ the vector space of traces on $\Ac$, while
the second associates to $\Ac$ its abelianization $\Ac/[\Ac, \Ac]$. 
From a conceptual viewpoint, it is important to realize 
the higher co/homologies $ HC^{*}€$, resp. $ HC_{*}$, as
derived functors. The obvious obstruction to such an 
interpretation is the non-additive nature of 
the category of algebras and algebra homomorphisms. 
This has been remedied in \cite{C3}, by
replacing it with the category of $\Lambda$-modules 
over the {\it cyclic category} $\Lambda$.
 
The cyclic category $\Lambda$ is a small category, obtained by
enriching with {\it cyclic morphisms} the familiar {\it simplicial
category} $\D$ of totally ordered finite sets and increasing maps. We
recall the presentation of $\D$ by generators and relations. 
It has one object $[n] = \{0 < 1 < \ldots < n \}$ for each
integer $n \geq 0$, and is generated by faces $\delta_i: [n-1] \ra [n]$ 
(the injection that misses $i$), and degeneracies $\s_j: [n+1] \ra [n] $
(the surjection which identifies $j$ with $j+1$), with the following 
relations:
\begin{equation}\label{ad}
\delta_j  \delta_i = \delta_i  \delta_{j-1}
\ \hbox{for} \ i < j  , \ \s_j  \s_i = 
\s_i  \s_{j+1} \qquad i \leq j 
\end{equation}
$$
\s_j  \delta_i = \left\{ \matrix{
\delta_i  \s_{j-1} \hfill &i < j \hfill \cr
1_n \hfill &\hbox{if} \ i=j \ \hbox{or} \ i = j+1 \cr
\delta_{i-1}  \s_j \hfill &i > j+1  . \hfill \cr
} \right.
$$
To obtain $\Lambda$ one adds for each $n$ a new morphism $\tau_n: [n]
\ra [n]$ such that,
\begin{equation}\label{ae}
\matrix{
\tau_n  \delta_i = \delta_{i-1}  
\tau_{n-1} &1 \leq i \leq n ,  \hfill \cr
\cr
\tau_n  \s_i = \s_{i-1} 
\tau_{n+1} &1 \leq i \leq n ,\cr
\cr
\tau_n^{n+1} = 1_n  . \hfill \cr
} 
\end{equation}
Note that the above relations also imply:
\begin{equation}\label{as}
 \tau_n  \delta_0 = \delta_n \, , \qquad \tau_n  \s_0 = \s_n  
 \tau_{n+1}^2 \, .
\end{equation}

Alternatively, $\Lambda$ can be defined by means of its ``cyclic
covering'', the category $E \Lambda$. The latter has one object $(\Zb ,
n)$ for each $n \geq 0$ and the morphisms $f : (\Zb , n) \ra (\Zb , m)$ are
given by non decreasing maps $f : \Zb \ra \Zb \ $, such that
$ f(x+n) = f(x)+m , \quad \fl  x \in \Zb$.
One has
$\Lambda = E  \Lambda / \Zb$,
with respect to the obvious action of $\Zb$
by translation.

To any algebra $A$ one associates a module $\Ac^{\natural}€$ over the category
$\Lambda$ by assigning to each integer $n \geq 0$ the vector space
$C^n (\Ac)$ of ($n+1$)-linear forms $\vp (a^0 , \ldots , a^n)$ on $A$, 
and to the generating morphisms the operators 
$\delta_i : C^{n-1} \ra C^n$, $\s_i : C^{n+1} \ra C^n$
defined as follows:
\begin{equation}\label{ag}
\matrix{
(\delta_i  \vp) (a^0 , \ldots , a^n) &=& \vp (a^0 , \ldots , a^i
a^{i+1} , \ldots , a^n), \quad i=0,1,\ldots , n-1 \, , \cr
\cr
(\delta_n  \vp) (a^0 , \ldots , a^n) &=& \vp (a^n  a^0 , a^1 , \ldots
, a^{n-1}) \, ; \hfill \cr
\cr
(\s_0  \vp) (a^0 , \ldots , a^n) &=& \vp (a^0 , 1 ,  a^1 , \ldots , a^n ) \, ,
 \hfill \cr
\cr 
(\s_j  \vp) (a^0 , \ldots , a^n) &=& \vp (a^0 , \ldots , a^j , 1, 
a^{j+1}€,
\ldots , a^n) , \quad j=1,\ldots , n-1 \, , \cr
\cr
(\s_n  \vp) (a^0 , \ldots , a^n) &=& \vp (a^0 , \ldots , a^n , 1)\, ; 
 \hfill \cr
\cr
(\tau_n  \vp) (a^0 , \ldots , a^n) &=& \vp (a^n , a^0 , \ldots ,
a^{n-1}) \, . \hfill \cr
}
\end{equation}
These operations satisfy the relations (\ref{ad}) and (\ref{ae}), which
shows that  $\Ac^{\natural}€$ is indeed a $\Lambda$-module. 

One thus obtains the desired interpretation of the cyclic
co/homology groups
of a $k$-algebra $\Ac$ over a ground ring $k$ in terms
of derived functors over the cyclic category (\cite{C3}): 
$$
HC^{n}€(\Ac) \simeq Ext_{\Lambda}^{n} (k^{\natural}€ , \Ac^{\natural}€) 
\quad {\rm and} \quad
HC_{n}€(\Ac) \simeq Tor^{\Lambda}_{n} (\Ac^{\natural}€ , k^{\natural}€) .
$$ 
Moreover, all of
the fundamental properties of the cyclic
co/homology of algebras, 
such as the long exact sequence relating it to Hochschild
co/homology (\cite{C2}, \cite{L}), are shared by the functors 
$Ext_{\Lambda}^{*}$/$Tor^{\Lambda}_{*}$-functors 
and, in this  generality, can be attributed to
the coincidence between the classifying space $B\Lambda$ of the
small category $\Lambda$ and the classifying space $BS^1 \simeq
P_{\infty}(\Cb)$ of the circle group. 

Let us finally mention that, from the very definition
of $Ext_{\Lambda}^{*} (k^{\natural}€ , F)$ and the existence of
a canonical
projective biresolution for $k^{\natural}$  (\cite{C3}), 
it follows that
the cyclic cohomology groups $HC^{*}€(F)$ 
of a $\Lambda$-module $F$, as well as the periodic ones
$HC_{\rm per}€ ^{*}€(F)$,
can be computed by means of a bicomplex analogous
to (\ref{abB}). A similar statement holds for the cyclic homology
groups.
\vspace{1cm}

\section{Cyclic theory for Hopf algebras}

The familiar antiequivalence between suitable categories of 
{\it spaces} and matching categories of {\it associative algebras},
effected by the passage to coordinates,
is of great
significance in both the purely algebraic context (affine schemes
versus commutative algebras) as well as
the topological one  (locally compact spaces versus commutative
$C^{*}€$-algebras). By extension, it has been adopted as a fundamental 
principle of noncommutative geometry. When applied to the realm of 
{\it symmetry}, it leads to promoting the notion of a {\it group}, 
whose coordinates form a commutative
Hopf algebra,
to that of a general {\it Hopf algebra}.
The cyclic categorical formulation recalled
above allows to adapt cyclic co/homology in a natural way 
to the treatment of symmetry in noncommutative geometry.
This has been done in \cite{CM2}, \cite{CM3} and will be
reviewed below.

We consider a Hopf algebra $\Hc$ over $k = \Rb \, \, {\rm or} \, \,  \Cb$,
with unit $\eta : k \ra \Hc$, counit $\ve : \Hc \ra k$ 
and antipode $S : \Hc \ra \Hc$. We use the standard definitions
(\cite{S}) together the usual convention for
denoting the coproduct:
\begin{equation}\label{bb}
\D(h)=  \sum  h_{(1)}  \ot  h_{(2)} \, , \qquad h \in \Hc .
\end{equation}
Although we work in the algebraic context, we shall include a datum
intended to play the r\^ole of the modular
function of a locally compact group. For reasons of consistency with 
the Hopf algebra context, this datum has a 
self-dual nature: it comprises
both a character $\delta \in {\Hc}^{*}€$,
\begin{equation}
\delta (a b) = \delta (a) \delta (b) \, , \quad \fl \, a, b \in 
\Hc ,
\end{equation}
and a group-like element $\sigma \in \Hc$,
\begin{equation}
\D(\sigma) = \s \ot \s \, , \qquad \ve (\s) = 1 ,
\end{equation}
related by the condition
\begin{equation}
    \delta (\s) = 1 .
\end{equation}    
Such a pair $(\delta, \sigma)$ will be called a {\it modular pair}.

\noindent The character $\delta$ gives rise to a  
$\delta$-{\it twisted antipode} $\wt S = S_{\d}€ : \Hc \ra \Hc $,
 defined by
\begin{equation}\label{dS}
\wt S (h) = \sum_{(h)} \delta (h_{(1)}) \ S (h_{(2)}) \quad , 
\quad h \in \Hc.
\end{equation}
Like the untwisted antipode, $\wt S$ is an 
algebra antihomomorphism
\begin{equation}
\matrix{
&\wt S (h^1  h^2) = \wt S (h^2)  \wt S (h^1) \quad , \quad \fl  
h^1 , h^2 \in \Hc \cr \cr
&\wt S (1) = 1, \hfill \cr
}
\end{equation}
a coalgebra twisted antimorphism
\begin{equation}
\D  \wt S (h) = \sum_{(h)} S (h_{(2)}) \ot \wt S (h_{(1)}) \quad , 
\quad \fl  h \in \Hc;
\end{equation}
and it also satisfies the identities
\begin{equation}\label{cd}
\ve \circ \wt S = \delta, \quad \quad 
\delta \circ \wt S = \ve.
\end{equation}

We start by associating to  $\Hc$, viewed only as a {\it coalgebra},   
the standard cosimplicial module known as the {\it cobar resolution}
({\cite{Ad}, \cite{Ca}), twisted by the insertion of the group-like
element $\s \in \Hc$. 
Specifically,  we set 
$C^n (\Hc) = \Hc^{\ot n}, \quad \fl n \geq 1 \quad {\rm and} 
\quad C^0 (\Hc) = k ,$
then define the {\it face} operators
$\delta_i: C^{n-1} (\Hc) \ra C^n (\Hc), \quad  0 \leq i \leq n,$
as follows: if ${n > 1}$, 
\begin{eqnarray}\label{ce}
&&\delta_0 (h^1 \ot \ldots \ot h^{n-1}) = 1 \ot h^1 
\ot \ldots \ot h^{n-1}, \nonumber \\
&& \nonumber \\
&&\delta_j (h^1 \ot \ldots \ot h^{n-1}) = h^1 \ot \ldots \ot \D h^j \ot 
\ldots \ot h^{n-1}  \nonumber \\
&&  = \sum_{(h_{j})} h^1 \ot \ldots \ot h_{(1)}^j \ot h_{(2)}^j \ot \dots \ot 
h^{n-1} ,
\quad \quad 1 \leq j \leq n-1,\nonumber\\
&& \\
&&\delta_n (h^1 \ot \ldots \ot h^{n-1}) = h^1 \ot \ldots \ot h^{n-1}
\ot \sigma , \nonumber
\end{eqnarray}
while if $n=1$
$$\delta_0 (1) = 1, \, \delta_1 (1) = \s .$$
Next, the {\it degeneracy} operators $\s_i : C^{n+1} (\Hc) \ra C^n (\Hc), 
\, \quad  0 \leq i \leq n,$ are defined by:
\begin{eqnarray}\label{cf}
&&\s_i (h^1 \ot \ldots \ot h^{n+1}) = h^1 \ot \ldots \ot \ve (h^{i+1}) 
\ot \ldots \ot h^{n+1}  \nonumber \\
&& \quad \quad = \ve (h^{i+1}) \, h^1 \ot \ldots \ot h^{i} \ot h^{i+2} 
\ot \ldots \ot h^{n+1}
\end{eqnarray}
and for $n=0$
$$\s_0 (h) = \ve (h), \, \, h \in \Hc . $$

\noindent The remaining features of the given data, namely
the {\it product} and the {\it antipode} of $\Hc$ together
with the character $\delta \in {\Hc}^{*}€$, are used to define the
candidate for the \textit{cyclic operator}, 
$\tau_n : C^n (\Hc) \ra C^n (\Hc) ,$ as follows:
\begin{eqnarray}\label{cg}
&& \tau_n (h^1 \ot \ldots \ot h^n) = (\D^{n-1}  \wt S (h^1)) \cdot h^2 
\ot \ldots \ot h^n \ot \sigma  \nonumber \\
&& \quad  = \sum_{(h^{1})} S(h^1_{(n)}) h^2 \ot \ldots \ot S(h^1_{(2)}) h^n 
\ot \wt S (h^1_{(1)}) \sigma .
\end{eqnarray}

\noindent Note that
$$\tau_{1}^{2}€ (h) =  \tau_{1} (\wt S (h) \sigma) = {\s}^{-1}€
{ \wt  S}^{2}€(h) \s ,$$
therefore the following is a necessary condition for
cyclicity:
\begin{equation}\label{ic}
(\sigma^{-1} \circ \wt{S})^2 = I.
\end{equation}
The remarkable fact is that this condition is also sufficient for
the implementation of the sought-for $\Lambda$-module.

A modular pair
${(\delta,\sigma)}$ satisfying (\ref{ic}) is called a  {\it modular pair
in involution}.

\begin{theorem}[\cite{CM2}, \cite{CM3}]
Let $\Hc$ be a Hopf algebra endowed
with a modular pair
${(\delta,\sigma)}$ in involution.
Then $\Hc_{(\delta,\sigma)}^{\natural} = \{ C^n (\Hc) \}_{n \geq 0}$ equipped 
with the operators given by (\ref{ce}) -- (\ref{cg}) is a module over the
cyclic category $\Lambda$.
\end{theorem}

The cyclic cohomology groups corresponding to the $\Lambda$-module  
$\Hc_{(\delta,\sigma)}^{\natural} \, $, denoted  
$H  C_{(\delta,\sigma)}^* (\Hc)$, 
can be computed from the bicomplex 
$( C C^{*, *} (\Hc), \, b , \, B ) ,$ analogous to (\ref{abB}),
defined as follows:
\begin{equation} \label{CbB} 
\matrix{ 
&C C^{p, q} (\Hc) = C^{q-p} (\Hc), \quad q \geq p , \cr \cr
& C C^{p, q} (\Hc) = 0 ,  \quad q < p \, ; \hfill \cr
}
\end{equation}
the operator 
\begin{equation} \label{Hb}
b: C^{n-1} (\Hc) \ra C^n (\Hc), \qquad    
b = \sum_{i=0}^{n} (-1)^i \delta_i \, .
\end{equation}
is explicitly given, if $n \geq 1$, by
\begin{eqnarray*}
\ b (h^1 \ot \ldots \ot h^{n-1})& = & 1 \ot h^1 \ot \ldots \ot h^{n-1}   \cr
& + & \ \sum_{j=1}^{n-1} (-1)^j
\sum_{(h_{j})} h^1 \ot \ldots \ot h_{(1)}^j \ot h_{(2)}^j \ot \dots \ot h^{n-1} \cr
&+ & \ (-1)^n h^1 \ot \ldots \ot h^{n-1} \ot \sigma  ,
\end{eqnarray*}
and if $n=0$ 
$$ b(c) = c \cdot (1 - \s ), \quad c \in k \, ;$$
the operator
$B: C^{n+1} (\Hc) \ra C^n (\Hc)$  
is defined by the formula
\begin{equation}\label{HB}
B =  N_n \circ \wt{\s}_{-1} \circ (1 + (-1)^n \tau_{n+1} ) , \quad n \geq 0 ,
\end{equation}
where $\wt{\s}_{-1} : C^{n+1} (\Hc) \ra C^n (\Hc) ,$
is the {\it extra degeneracy} operator 
\begin{eqnarray}\label{ex}
&& \wt{\s}_{-1} (h^1 \ot \ldots \ot h^{n+1}) = (\D^{n-1}  \wt S (h^1)) \cdot h^2 
\ot \ldots \ot h^{n+1} \nonumber \\
&& \quad = \sum_{(h^{1})} S(h^1_{(n)}) h^2 \ot \ldots \ot S(h^1_{(2)}) h^n 
\ot \wt S (h^1_{(1)})  h^{n+1} ,
\end{eqnarray}
$$ \wt{\s}_1 (h) = \d (h), \, \, h \in \Hc $$
and
\begin{equation}\label{N}
N_n = 1 + (-1)^n \tau_n  + \ldots + (-1)^{n^{2}}€ {\tau_n}^n.
\end{equation}
Explicitly,
\begin{eqnarray}\label{norm}
N_n (h^1 \ot \ldots \ot h^n) =  \quad \quad \quad \quad \quad \quad \quad \quad  
\quad \quad \quad \quad \\
\sum_{j=0}^{n} (-1)^{nj} \D^{n-1} \wt S (h^j) \cdot h^{j+1} 
\ot \ldots \ot h^n \ot \s \ot \wt S^2 (h^{0}) \s \ot \ldots \ot \wt S^2 (h^{j-1}) \s
\nonumber \\
=\sum_{j=0}^{n} (-1)^{nj} \D^{n-1} \wt S (h^j) \cdot h^{j+1} 
\ot \ldots \ot h^n \ot \sigma \ot \s h^{0} \ot \ldots \ot \s h^{j-1} \quad \quad \quad
\nonumber \\
= \sum_{j=0}^{n} (-1)^{nj} \sum_{(h_{j})} S (h^j_{(n)}) \ot \ldots \ot 
S (h^j_{(2)}) \ot \wt S (h^j_{(1)}) \cdot \quad \quad \quad \quad \quad \quad \quad
\quad \nonumber \\
\quad \cdot h^{j+1} \ot \ldots \ot h^n \ot \s \ot \s h^0 \ot \ldots \ot 
\s h^{j-1} . \nonumber
\end{eqnarray}
\bigskip
In particular, for $n=0$,
\begin{equation}\label{B0}
B(h) = \delta (h) + \ve (h) , \quad \fl \ h \in \Hc .
\end{equation}
The expression of the $B$-operator can be simplified by passing to the 
quasi-isomorphic \textit{normalized bicomplex}
$( C  {\bar C}_{(\delta,\sigma)}^{*, *} (\Hc), \, b , \, {\bar B}) ,$ 
defined as follows
\begin{equation}\label{Bbar}
\matrix{ 
&C  {\bar C}_{(\delta,\sigma)}^{p, q} (\Hc) = {\bar C}^{q-p} (\Hc), \quad q \geq p , 
\cr \cr
& C  {\bar C}_{(\delta,\sigma)}^{p, q} (\Hc) = 0 ,  \quad q < p \, ,\hfill \cr
}
\end{equation}
where
$$ {\bar C}^n (\Hc) =  ({\rm Ker} \ve)^{\ot n}, \quad \fl n \geq 1, \quad \quad 
			{\bar C}^0 (\Hc) = k ;$$
while the formula for the $b$-operator remains unchanged,
the new horizontal operator becomes 
\begin{equation}\label{barB}
\bar B =  N_n \circ \wt{\s}_{-1} , \quad \quad n \geq 0 ;
\end{equation}
in particular, for $n=0$, one has
\begin{equation}\label{nB0}
{\bar B}(h) = \delta (h) , \quad \fl \ h \in \Hc .
\end{equation}

\noindent An alternate description of the cyclic cohomology
groups $H  C_{(\delta, 1)}^* (\Hc)$, in terms of the
Cuntz-Quillen formalism, is given by M. Crainic in \cite{Cr}.
\medskip

We should also mention that the corresponding
{\it cyclic homology} groups
\begin{equation} \label{CH}
H  C^{(\delta, \s)}_{n}€ (\Hc) =
 Tor^{\Lambda}_{n} (\Hc_{(\delta,\sigma)}^{\natural}€ , k^{\natural}€)
\end{equation}
can be computed from the
bicomplex $( C C_{*, *} (\Hc), \, b , \, B ) ,$ obtained
by dualising
(\ref{CbB}) in the obvious fashion:
\begin{equation}\label{HbB}
\matrix{ 
&C C_{p, q} (\Hc) = Hom (C^{q-p} (\Hc), k) , \quad q \geq p , \cr \cr
& C C_{p, q} (\Hc) = 0 ,  \quad q < p \, , \hfill \cr
}
\end{equation}
with the boundary operators $b$ and $B$ the transposed of 
the corresponding coboundaries.
\medskip

When applied to the usual notion of symmetry in differential
geometry, the ``Hopf algebraic'' version of cyclic cohomology 
discussed above recovers both the Lie algebra co/homology and
the differentiable cohomology of Lie groups, as illustrated by
the following results.

\begin{proposition}[\cite{CM2}] Let $\FG$ be a Lie algebra and
let   $\delta : \FG \ra \Cb$ be a character of $\FG$. 
With $\Uc (\FG) $ denoting
the enveloping algebra of ${\FG}_{\Cb}€$, viewed
as a Hopf algebra with modular pair $(\delta, 1)$, one has 
$$
 H  C_{{\rm per} \, (\delta, 1)}^* \, (\Uc (\FG)) \simeq 
 \sum^{\oplus}€_{i \equiv * \, (2)}€ \quad H_i \, ({\FG} , {\Cb}_{\delta}€) ,
$$
where ${\Cb}_{\delta}€$ is the $1$-dimensional $\FG$-module associated
to the character $\delta$.
\end{proposition}

\begin{remark} {\rm In a dual fashion, one can prove that}
 \begin{equation}  \label{DL}
 H  C^{{\rm per} \, (\delta, 1)}_{*}€ \, (\Uc (\FG)) \simeq 
 \sum^{\oplus}€_{i \equiv * \, (2)}€ \quad 
 H^{i} \, ({\FG} , {\Cb}_{\delta}€) .
 \end{equation}
\end{remark}

\begin{proposition}[\cite{CM2}] Let ${\Hc} (G)$ be 
the Hopf algebra of polynomial functions
on a simply connected affine algebraic nilpotent group $G$,
with Lie algebra $\FG$. 
Then its periodic cyclic cohomology with respect to the trivial
modular pair $(\ve, 1)$ coincides with the Lie algebra cohomology:
$$
 H  C_{\rm per}^* \, ({\Hc} (G)) \simeq  
 \sum^{\oplus}€_{i \equiv * \, (2)}€ \quad H^{i} \, ({\FG}, \Cb) . 
$$ 
\end{proposition}

\begin{remark} {\rm Since by Van Est's Theorem, the
cohomology of the nilpotent Lie algebra
$\FG$ is isomorphic to the differentiable group cohomology
$H_{\rm d}^* (G)$, the above isomorphism can be reformulated as}
$$
 H  C_{\rm per}^* \, ({\Hc} (G))) \simeq  
 \sum^{\oplus}€_{i \equiv * \, (2)}€ \quad H_{\rm d}€^{i} \, (\FG, \Cb) . 
$$
Under the latter form it continues to hold for any affine
algebraic Lie group $G$, {\rm  with the same proof as in \cite{CM2} }.
\end{remark}
\smallskip

Hopf algebras often arise implemented as endomorphisms
of associative algebras. A {\it Hopf action} of a Hopf algebra $\Hc$ on
an algebra $\Ac$ is given by a linear map,
$\Hc \ot \Ac \ra \Ac, \quad h \ot a \ra h(a) $
satisfying the {\it action property}
\begin{equation} \label{AP}
h_1 (h_2  a) = (h_1  h_2) (a), \qquad \fl  h_i \in {\Hc}, \,
a \in {\Ac}
\end{equation}
and the {\it Hopf-Leibniz rule}
\begin{equation} \label{HL}
h(ab) = \sum  h_{(1)}  (a) \, h_{(2)}  (b),  \qquad \fl  a,b  \in {\Ac},
\, h \in {\Hc} .
\end{equation}

In practice, $\Hc$ first appears as a subalgebra
of endomorphisms of an algebra
$\Ac$ fulfilling  (\ref{HL}), and it is precisely the Hopf-Leibniz rule 
which dictates the coproduct of $\Hc$. In turn, the modular pair 
of $\Hc$ arises in connection with the existence
of a twisted trace on $\Ac$.

Given a Hopf action $\Hc \ot \Ac \ra \Ac $ together with a modular
pair $(\delta, \s)$, a
linear form $\tau : \Ac \ra \Cb$ is called a $\sigma$--{\it trace}
under the action of $\Hc$ if 
\begin{equation} \label{tt}
\tau (ab) = \tau (b \sigma (a)) \qquad \fl  a,b  \in \Ac  .
\end{equation}
The $\sigma$--trace $\tau$ is called $\delta$--{\it invariant}
under the action of $\Hc$ if 
\begin{equation}\label{ip}
\tau (h(a)b) = \tau (a  \wt S (h)(b)) \qquad \fl  a,b  \in \Ac  , \ h
\in {\Hc}. 
\end{equation}
If $\Ac$ is unital, the ``integration by parts'' formula (\ref{ip}) is
equivalent to the $\delta$--invariance condition
$$
\tau (h(a)) = \delta (h) \, \tau (a)
\qquad \fl  a  \in \Ac  , \ h \in {\Hc}. 
$$

With the above assumptions, the very definition of the cyclic
co/homology of $\Hc$, with respect to a modular pair in involution
$(\delta, \s)$, is uniquely dictated such that the following
{\it Hopf action principle} holds:

\begin{theorem}[ \cite{CM2},  \cite{CM3}] 
    Let $\tau : \Ac \ra \Cb$ be a $\delta$--invariant 
 $\sigma$--trace under the Hopf action of $\Hc$ on $\Ac$. Then 
the assignment 
\begin{eqnarray} \label{cm}
  \g (h^1 \ot \ldots \ot h^n) (a^0 , \ldots , a^n) = 
 \tau (a^0  h^1 (a^1) \ldots h^n (a^n)) ,  \\
 \fl \quad  a^0 , \ldots , a^n \in \Ac, \qquad
 h^1 , \ldots , h^n \in \Hc \nonumber
\end{eqnarray}
defines a map of $\Lambda$--modules 
$ \g^{\natural}€: \Hc_{(\delta,\sigma)}^{\natural} \ra  {\Ac}^{\natural}€$,
which in turn induces characteristic
 homomorphisms in cyclic co/homology:
 \begin{eqnarray} 
     \g_{\tau}^{*}€ : H  C_{(\delta,\sigma)}^* (\Hc) \ra
     H  C^* (\Ac) ; \label{cco} \\
     \g^{\tau}_{*}€€ :  H  C_{*}€ (\Ac) \ra H C^{(\delta, \s)}_{*}€ (\Hc) .
     \label{cho}
 \end{eqnarray}
\end{theorem}
\smallskip

As a quick illustration,  
let us assume that $\FG$ is a Lie algebra of derivations
of an algebra $\Ac$ and $\tau$ is a  $\delta$--invariant 
trace on $\Ac$. 
Then (\ref{cho}), composed with
the Chern character in cyclic homology (\ref{cch}) on one hand
and with the isomorphism (\ref{DL}) on the other, 
recovers the additive map 
\begin{equation} \label{Kco}
    ch_{\tau}€^{*}€ : K_{*}€(\Ac) \ra H^{*}€ ({\FG} , {\Cb}_{\delta}€) ,
\end{equation}
previously introduced in \cite{C1}, in terms of
$\FG$--invariant curvature forms associated to 
an arbitrary $\FG$--connection.
Before applying the isomorphism (\ref{DL}), the periodic cyclic class
$  \g_{*} (ch_{*} (e) ) \in 
H  C^{{\rm per} \, (\delta, 1)}_{*}€ \, (\Uc (\FG)), \, $ for
$e^{2}€ = e \in \Ac $, 
is given by the cycle with the following components:
\begin{eqnarray} \label{eco}
  \g_{*} (ch_{0} (e) )  =  \tau (e) \qquad {\rm and}  \qquad \fl k \geq 1,
 \qquad \fl  h^{1}€, \ldots ,  h^{2k}€ \in \Uc (\FG)  \nonumber \\
   \g_{*} (ch_{2k} (e) )  \, (h^{1}€, \ldots , h^{2k}€) = 
  \quad \quad \quad \quad \quad \quad \quad \\
 (-1)^{k}€{(2k)! \over k!}  \left( \tau (e h^{1}€(e) \cdots h^{2k}€(e) ) -
  {1 \over 2} \tau ( h^{1}€(e) \cdots h^{2k}€(e) ) \right) . \nonumber
\end{eqnarray} 
The Lie algebra cocycle representing the class
$ ch_{\tau}€^{*}€(e) \in  H^{*}€ ({\FG} , {\Cb}_{\delta}€)$
in terms of the Grassmannian connection
is obtained by restricting $ \g^{\tau}_{*} (ch_{*} (e) )$
to $\wedge^{\cdot}€ \FG$ via antisymmetrization.

\vspace{1cm}

\section{Transverse index theory on general foliations}

The developments discussed in the preceding section have been 
largely motivated by a challenging computational problem concerning 
the index of transversely hypoelliptic differential operators
on foliations \cite{CM1}. In turn, they were instrumental in 
settling it \cite{CM2}. The goal of this section is to
highlight the main steps involved.
 
The transverse geometry of a foliation $(V, \Fc)$, i.e. the geometry of
the ``space'' of leaves $V/\Fc$, provides a prototypical example of 
noncommutative space, which already exhibits many of the distinctive
features of the general theory.  
In noncommutative geometry, a  geometric space is given by a {\it
spectral triple} $(\Ac ,\Hc ,D)$,
where $\Ac$ is an involutive algebra of operators in a Hilbert space
$\Hc$, representing the ``local coordinates'' of the space,
and $D$ is an unbounded selfadjoint operator on $\Hc$.  
The operator $D^{-1} = ds$
corresponds to the infinitesimal line  element in Riemannian geometry 
and, in addition to
its metric significance, it also carries nontrivial homological
meaning, representing a $K$-{\it homology} class of $\Ac$.
The construction of such a spectral triple associated to a general 
foliation  (\cite{CM1}) comprises several steps and incorporates 
important ideas from \cite{C4} and \cite{HS}. 

To begin with, we recall that a codimension $n$
foliation $\Fc$ of an $N$-dimensional manifold $V$
can be given by means of a {\it defining cocycle}
$( U_{i}€ , f_{i} , g_{ij} )$, where
$\{U_{i}€\}$ is an open cover of $V$, $f_{i}€: U_{i}€ \ra T_{i}€$
are submersions with connected fibers onto 
$n$-dimensional manifolds $\{T_{i}€\}$ and 
$$
g_{ij}€ : f_{j}€( U_{i} \cap U_{j}) \ra f_{i}€( U_{i} \cap U_{j}) 
$$
are diffeomorphisms such that
$$
   \fl \quad (i, j) \, , \qquad
    f_{i}€ = g_{ij} \circ f_{j}€  \qquad {\rm on} \quad
    U_{i} \cap U_{j} \, .
$$
The disjoint union $M = \cup_{i}€ \, U_{i}€ \times \{i \}$ 
can be regarded as a  complete transversal
for the foliation, while the collection of local diffeomorphisms
$\{ g_{ij} \}$ of $M$
generates the transverse {\it holonomy pseudogroup} $\Gamma$. 
We should note that the notion of {\it  pseudogroup} used here 
is slightly different from the standard one, since we do not 
enforce the customary hereditary condition; in particular, any
group of diffeomorphisms is such a pseudogroup.

We shall assume $\Fc$ {\it transversely oriented}, which
amounts to stipulating that $M$ is oriented and that
$\Gamma$ consists of orientation preserving local diffeomorphisms.
From $M$ we shall pass by a ${\rm Diff}^{+}€$--functorial 
construction (\cite{C4})
to a quotient bundle, 
$\pi: P \ra M $, of the frame bundle, whose sections are the
Riemannian metrics on $M$. 
Specifically, $P = F /SO (n) $, where $F$ is
the $GL^{+}€(n, \Rb)$--principal bundle of oriented frames on
$M$. The total space $P$ admits a canonical {\it para-Riemannian
structure} as follows. 
The vertical subbundle $\Vc \sbs TP$, $\Vc =\Ker \pi_*$, 
carries natural Euclidean structures on each of its
fibers, determined solely by the choice of a $GL^+(n,  \Rb)$--invariant
Riemannian metric on 
the symmetric space $GL^{+}€(n,  \Rb) / SO (n)$. On the other hand,
the quotient bundle $\Nc = (TP)/ \Vc$ comes equipped with a
tautologically defined Riemannian structure: every $p\in P$ is an
Euclidean structure on $T_{\pi (p)} (M)$ which is identified to 
$\Nc_{p}€$ via $\pi_*$. 
\smallskip

The naturality of the above construction with respect to 
${\rm Diff}^{+}€$
ensures that the action of the holonomy pseudogroup $\Gamma$
lifts to both bundles $F$ and $P$. One can thus form for each the 
associated {\it smooth \'etale groupoid} 
$ F \semi \G$,  resp. $ P \semi \G $ .
An element of $F \semi \G$ or of $ P \semi \G $ is given
by a pair
$$ (x,\vp) \, , \qquad x \in {\rm Range} \, \vp \, ,
$$ 
while the composition law is
$$
(x,\vp) \circ (y,\psi) = (x,\vp \circ \psi) \quad \hbox{if} \quad 
y \in {\rm Dom} \, \vp \quad \hbox{and} \quad \vp (y) = x \, . 
$$
We let 
$$ {\Ac}_{F}€ = C_c^{\ify} ( F \semi \G) \, , \quad {\rm resp.} \quad
     {\Ac}_{P}€ = C_c^{\ify} ( P \semi \G)   
$$
denote the corresponding convolution algebras. 
The elements of $ \Ac = {\Ac}_{P}€ $ will 
serve as ``functions
of local coordinates'' for the noncommutative space
$V/\Fc$. A generic element of $\Ac$ can be represented as a
linear combination of monomials
$$  a  = 
f \, U_{\psi}^* \, , \ f \in C_c^{\ify} ({\rm Dom} \, \psi) \, ,
$$
where the star indicates a contravariant notation. 
The multiplication rule is
$$
f_1 \, U_{\psi_1}^* \cdot f_2 \, U_{\psi_2}^* = 
f_1 \cdot  (f_2 \circ \wt{\psi}_1) \, U_{\psi_2 \psi_1}^* \, ,
$$
where by hypothesis the support of $f_1 (f_2 \circ \wt{\psi}_1)$ is a
compact subset of
$$
{\rm Dom} \, \psi_1 \cap \psi_1^{-1} \, {\rm Dom} \, \psi_2 \sbs {\rm Dom}
\,
\psi_2 \, \psi_1 \, . 
$$
The algebras $\Ac = {\Ac}_{P}€$ and ${\Ac}_{F}€$ 
admit canonical ${\ast}$--representations
on the Hilbert spaces
$$  L^{2}€(P) = L^{2}€(P, vol_{P}€) \, , \quad {\rm resp.} \quad
     L^{2}€(F) = L^{2}€(F, vol_{F}€) \, ,
$$
where $vol_{P}$,     
resp. $vol_F$ denotes the canonical ${\rm Diff}^{+}€$--invariant
volume form on $P$, resp. on $F$.
Explicitly, 
for $\Ac = {\Ac}_P $,
\begin{equation}
    ( (f \, U_{\psi}^*) \, \xi) (p) = f(p) \, \xi (\psi (p)) \quad \fl
\, p
\in P \, , \ \xi \in L^{2}€(P) \, , \ f \, U_{\psi}^* 
\in {\Ac} \, , 
\end{equation}
and similarly for ${\Ac}_{F}€$.
We shall denote by $A = \bar{\Ac}$ the norm closure of $\Ac$ in this 
representation.
\smallskip

Evidently, the algebra $\Ac$ depends on the choice of the defining
cocycle $( U_{i}€ , f_{i} , g_{ij} )$. However, if 
$( U^{\prime}€_{i}€ , f^{\prime}€_{i} , g^{\prime}€_{ij} )$
is another cocycle defining the same foliation  $\Fc$, 
with corresponding
algebra ${\Ac}^{\prime}€$ (resp.  $A^{\prime}€$ ), then $\Ac$ and
${\Ac}^{\prime}€$ are {\it \, Morita equivalent}, while
the $C^{*}€$--algebras 
$A$ and $A^{\prime}€$ are {\it \, strongly Morita equivalent}.
We recall that Morita equivalence preserves the cyclic co/homology
and the $K$-theory/$K$-homology. Also, in the commutative case
it simply reduces to isomorphism; e.g., absent any nontrivial
pseudogroup of diffeomorphisms,
two manifolds
$N$ and $N^{\prime}€$ are diffeomorphic iff the algebras
$C_c^{\ify} (N)$ and $C_c^{\ify} (N^{\prime}€)$ are Morita
equivalent.
\smallskip

To complete the description of the spectral triple associated to
$V/\Fc$, we need to define the operator $D$. In practice,
it is more convenient to work with another representative
of the same $K$-homology class, namely
the {\it hypoelliptic signature operator} $Q = D |D|$. The latter
is a second order differential operator, acting on the Hilbert space
\begin{equation}
 \Hc = {\Hc}_{P}€ : = L^{2}€({\wedge}^{\cdot}€ \Vc^{*}€  \ot 
 {\wedge}^{\cdot}€ \Nc^{*}€  , \,  vol_{P}€) \, ; 
\end{equation}
it is defined as a graded sum
\begin{equation} \label{Q}
    Q = (d_V^* \, d_V - d_V \, d_V^*) \op (d_H + d_H^*) \, ,
\end{equation}
where $d_V$ denotes the vertical exterior derivative and
$d_H $ stands for the horizontal exterior differentiation with
respect to a fixed connection on the frame bundle.
When $n \equiv 1 \, {\rm or} \, 2 \, ({\rm mod} \,  4)$, 
for the vertical component to make sense,
one has to replace $P$ with $ P \times S^1$ 
so that the dimension of the vertical fiber
stays even.

{\begin{proposition}[\cite{CM1}]. For any $a \in \Ac$, 
$[D,a]$ is bounded. 
For any $f \in C_{c}€^{\infty}€(P)$ and $ \lb \notin \Rb $,  
$\, f(D-\lb)^{-1}$ is $p$-summable, 
$\fl \, p > m = {n(n+1) \over 2} + 2n$.
\end{proposition}

One now confronts a well posed index problem. The operator $D$ 
determines an {\it index pairing} map
$ {\rm Index}_{D}€: K_{*}€ (\Ac) \ra \Zb $, as follows:
\begin{itemize}
    \item[(0)] in the {\it graded} (or {\it even}) case,
    $$ {\rm Index}_{D}€ ([e]) = {\rm Index} \, (e D^{+}€e) \,  , 
    \quad e^{2}€ = e \in \Ac \, ;
    $$
    \item[(1)] in the {\it ungraded} (or {\it odd}) case,
    $$ {\rm Index}_{D}€ ([u]) = {\rm Index}\, (P^{+}€ u P^{+}€) \,  , \quad
    u \in GL_{1}€(\Ac) \, ,
    $$
    where $P^{+}€ = {1+F \over 2}$, with  $F = \Sign \, (D)$ .
\end{itemize} 
One of the main functions of cyclic co/homology, its {\it raison d'\^etre}
in some sense, is to compute the index pairing via the equality
\begin{equation} \label{index}
 {\rm Index}_{D}€ (\kappa) = \lgl {ch}_* (D) , \, {\rm ch}^* (\kappa) \rgl 
\qquad \fl \, \kappa \in  K_{*}€ (\Ac)  .
\end{equation}
The cyclic cohomology class ${ch}_* (D) \in H C_{\rm per}^{*} (\Ac)$,
i.e. its {\it Chern character in $K$-homology}, is defined in 
the {\it ungraded case} by means of the cyclic cocycle
\begin{equation}  \label{Kch}
\tau_{F}€ (a^0 ,\ldots ,a^n) = \Trace \, (a^0 [F,a^1] \ldots [F,a^n]) 
\, , \qquad  a^j \in \Ac 
\end{equation}
where $n$ is any odd integer
exceeding the dimension of the spectral triple  $(\Ac ,\Hc , D)$;
in the {\it graded case},
$\Trace$ is replaced with the graded trace $\Trace_{s}€$ and $n$ is even.
Being defined by means of the {\it operator trace}, 
the cocycle (\ref{Kch}) 
is inherently difficult to compute. The problem is therefore
to {\it provide an explicit formula 
for the Chern character in $K$-homology}. 

We should note at this point that, for smooth groupoids
such as those associated to foliations,
the answer to (\ref{index})
is indeed known for all {\it $K$--theory classes in the range of the
assembly map} from the corresponding geometric $K$-group to 
the analytic one (cf. \cite{Co}).

As mentioned before, the functors $K$--theory/$K$--homology and 
cyclic co/homology are Morita invariant. 
Moreover, the corresponding Chern characters are {\it Morita
equivariant}, in such a way that both sides of (\ref{index}) 
are preserved
by the canonical isomorphisms associated with a
Morita equivalence datum.
Thus, one may as well take advantage of the Morita invariant nature of 
the problem and choose from the start 
a defining cocycle $( U_{i}€ , f_{i} , g_{ij} )$ for $\Fc$
with all local transversals $T_{i}€$  {\it flat affine manifolds}.
This renders $M$ itself as a flat affine manifold, although it
{\it does not} allow one to assume that the affine structure is
preserved by $\Gamma$. One can however take the horizontal
component in (\ref{Q}) with respect to a {\it flat connection} $\nb$. 
It is then readily seen that the operator $Q$ belongs to the
class of operators of the form
\begin{equation} \label{Ru}
    R = \pi_{a}€ (R_{\Uc}), \quad {\rm with} \quad
   R_{\Uc}€ \in (\Uc (G_{a}€(n)) \ot {\rm End} (E) )^{SO(n)}  \, ,
\end{equation}
where $ G_{a}€(n) = \Rb^{n}€ \semi GL(n, \Rb)$ is
the affine group, $\pi_{a}€$ denotes its
right regular representation and $E$ is a unitary 
$SO(n)$--module. 

A differential operator $R$ of the form (\ref{Ru}) will be called 
{\it affine}.
If in addition the principal symbol of $R$, with respect to 
(\ref{sym}),
is invertible then $R$ will be called an {\it hypoelliptic
affine operator}.

By an easy adaptation of a classical theorem of Nelson and Stinespring, 
one can show that {\it any  hypoelliptic
affine operator $R$ which is formally selfadjoint
is in fact essentially selfadjoint}, with core any dense
$ G_{a}€(n)$--invariant subspace of
$C^{\ify}$--vectors for $\pi_{a}€$.
\smallskip

The hypoelliptic calculus adapted to the para-Riemannian structure 
of the manifold $P$ and to the treatment of the above operators 
is a
particular case of the  pseudodifferential calculus 
on Heisenberg manifolds (\cite{BG}).  
One simply modifies the notion of  homogeneity of symbols
$\s (p,\xi)$ by using  the homotheties
\begin{equation}  \label{sym}
\lb \cdot \xi = (\lb \xi_v , \lb^2 \xi_n) \ , \ \fl \, \lb \in \Rb_+^* ,
\end{equation}
where $\xi_v$, $\xi_n$ are the vertical, resp. normal
components of
the covector $\xi$. The above formula depends on local coordinates 
$(x_v, x_n)$
adapted  to the vertical foliation,
but the corresponding pseudodifferential
calculus is independent of such choices. The principal symbol of a
hypoelliptic operator of order $k$ is a function 
on the fibers of ${\Vc}^* \op {\Nc}^*$,
homogeneous of degree $k$ in
the sense of (\ref{sym}). 
The distributional kernel
$k(x,y)$
of a pseudodifferential operator $T$ in this hypoelliptic calculus has  the
following behavior near the diagonal:
\begin{equation} \label{kxy}
k(x,y) \sim \sum a_j (x,x-y) - a(x) \log \vert x-y\vert' + O(1) \, ,
\end{equation}
where $a_j$ is homogeneous of degree  $\, -j$ in $\, x-y$ in the sense of
(\ref{sym}),
and the metric $\vert x-y\vert'$  is locally of the form
\begin{equation} \label{mxy}
\vert x-y \vert' = ((x_v - y_v)^4 + (x_n - y_n)^2)^{1/4}  \, .
\end{equation}
The 1--density $a(x)$ is independent on the choice of metric 
$\vert \cdot \vert'$ and can be obtained from the symbol of
order  $\, -m$ of $T$, 
where 
$$m= {n(n+1) \over 2} + 2n
$$ 
is the {\it Hausdorff
dimension} of the metric space $(P,\vert \cdot \vert')$.
Like in the ordinary pseudodifferential calculus, this allows to
define a residue of Wodzicki-Guillemin-Manin type, extending the
Dixmier trace to operators of all degrees, by the equality
\begin{equation} \label{res}
{\int \!\!\!\!\!\! -} T = {1 \over m} \int_{P}€ a(x) \, .
\end{equation}
\smallskip

One uses the hypoelliptic calculus to prove (\cite{CM1}) that
the spectral triple $(\Ac, \Hc, D)$, or more generally
that obtained by replacing $Q$
with any hypoelliptic affine operator $R$ (in which case $D |D| = R$),
fulfills the hypotheses of the
operator theoretic {\it local index theorem} of \cite{CM1}. 
Its application
allows to express the corresponding Chern character
$$ {ch}_* (R) = {ch}_* (D) \in H C_{\rm per}^{*} (\Ac_{P}€)
$$
in terms of
the {\it locally computable} residue (\ref{res}). In the odd case,
it is given by the cocycle $ \Phi_{R}€ = \{ \vp_n  \}_{n=1,3,\ldots} \,$ 
in the
$(b,B)$--bicomplex of $\Ac$ defined as follows:
\begin{eqnarray} \label{phi}
    \matrix{
\vp_n (a^0 ,\ldots ,a^n) =   \cr
\cr
 \sum_k c_{n,k} \,
{\int \!\!\!\!\!\! -}  a^0 [R, a^1]^{(k_1)} \ldots
[R, a^n]^{(k_n)}  \, \vert R \vert^{-n-2\vert k\vert} \, ,
\quad a^j \in \Ac , \cr
}
\end{eqnarray}
where we used the abbreviations
$$T^{(k)} = \nb^k (T) \qquad {\rm and} \qquad \nb (T) = D^2 T
- TD^2 \, ,
$$
$k$ is a multi-index, $\vert k \vert = k_1 +\ldots + k_n$,
and 
$$
c_{n,k} =
(-1)^{\vert k \vert} \, \sqrt{2i} (k_1 ! \ldots k_n
!)^{-1} \, ((k_1 +1) \ldots (k_1 + \ldots + k_n
+n))^{-1} \, \G ( \vert k \vert + {n\over 2} ) ;
$$
there are finitely many nonzero terms in the above sum and only
finitely many components of $\Phi_{R}€$ are nonzero.
In the even case, the corresponding cocycle 
 $ \Phi_{R}€ = \{ \vp_n  \}_{n=0,2,\ldots} \,$ is defined in
a similar fashion, except for $\vp_0$  (see \cite{CM1}).
\smallskip

The expression (\ref{phi}) is definitely  explicitly computable,
but its actual computation is exceedingly difficult to perform.
Already in the case of codimension $1$ foliations, where we did 
carry through its calculation by hand, it involves computing
thousands of terms. On the other hand, in the absence
of a guiding principle, computer calculations are unlikely to produce
an illuminating answer. 
However, a simple inspection of (\ref{phi}) reveals some helpful
general features. For the clarity of the exposition, we shall 
restrict our comments to the case $R =Q$, which is our  main case 
of interest anyway.

First of all, since the passage from $(\Ac_{F}€, \Hc_{F}€)$ to
$(\Ac_{P}€, \Hc_{P}€)$ involves the rather harmless operation
of taking $K$--invariants with respect to the compact group
$K = SO(n)$, 
we may work directly at the
level of the frame bundle, equivariantly with respect to $K$.
Secondly, since we are interested in the flat case,
we may as well assume for starters that
$M = \Rb^{n}€$, with the trivial
connection. This being the case,
we shall identify $F$ with the affine group $G_{a}€(n)$.
We may also replace $\G$ by the full group 
${\rm Diff}^{+}€ (\Rb^{n}€)$ and thus work for awhile with the algebra
$$\Ac (n) = C_{c}^{\infty}€€ (F(\Rb^{n}€) 
\semi {\rm Diff}^{+}€(\Rb^{n}).
$$

We recall that $Q$ is built from the vertical
vector fields $\{ Y_{i}€^{j}€ ; \, i, j = 1, \ldots , n \}$
which form the canonical basis of ${\frak gl} (n, \Rb)$ and
the horizontal vector fields $\{ X_{k}€ , \, k = 1, \ldots , n \}$
coming from the canonical basis of $\Rb^{n}€$. Therefore, 
the expression under the residue-integral in (\ref{phi}) involve
iterated commutators of these vector fields with multiplication
operators of the form $ a  = f \, U_{\psi}^* \, , 
\ f \in C_c^{\ify} (F, ) \, \psi \in \Gamma $.
Now the canonical action of $GL^+ (n, \Rb)$ on $F$ commutes 
with the action
of $\G$ and hence extends canonically to the crossed product
$\Ac_{F}€$. At the Lie
algebra level, this implies that the operators on $\Ac_{F}€$
defined by
\begin{equation} \label{Y}
Y_{i}^j (f \, U_{\psi}^*) = (Y_i^{j} \, f) \, U_{\psi}^* \, 
\end{equation}
are derivations:
\begin{equation} \label{YR}
Y_{i}^j (ab) = Y_i^{j} (a) \, b + a \, Y_{i}^j (b) \, .
\end{equation}
The horizontal
vector fields $X_k$ on $F$ can also be made to act on the crossed
product algebra, according to the rule
\begin{equation} \label{X}
X_k (f \, U_{\psi}^*) = X_k (f) \, U_{\psi}^* \, . 
\end{equation}
However, since the trivial connection is not preserved by
the action of $\G$, the operators  $X_k$ are no longer 
derivations of  $\Ac_{F}€$;
they satisfy instead
\begin{equation} \label{XR}
X_i (ab) = X_i (a) \, b + a \, X_i (b) + \sum \d_{ji}^k (a) \, Y_k^j (b) \, ,
\, 
a,b \in {\Ac} \, .
\end{equation}
The linear operations $\d_{ij}^k$ are of the form
\begin{equation} \label{dij}
 \d_{ij}^k (f \, U_{\psi}^*) = \g_{ij}^k \, f \, U_{\psi}^* \, ,
\end{equation}
with $\g_{ij}^k \in  C^{\infty}€ (F \semi \G) $ 
characterized by the identity
\begin{equation} \label{gij}
{\psi}^* \, \om_{j}^{i}€ - \om_{j}^{i}€ =  \sum_{k}€ \g_{jk}^{i}  \, 
\theta^k \, ,
\end{equation}
where $ \om $ is the standard {\it flat connection form} 
and $\theta $ is the {\it fundamental form} on $F = F(\Rb^{n}€)$.
From (\ref{gij}) it easily follows that each  $\d_{ij}^k$ 
is a derivation:
\begin{equation} \label{dijR}
 \d_{ij}^k (ab) = \d_{ij}^k (a) \, b + a \, \d_{ij}^k (b) \, . 
\end{equation}
The commutation of the $Y_i^{j}$ with $\d_{ab}^c$ 
preserves the linear span of the $\d_{ab}^c$.
However, the successive
commutators with $X_k$ produces new operators
$$
\d_{ab,i_1 \ldots i_r}^c = [X_{i_r} , \ldots [X_{i_1} , \d_{ab}^c] \ldots ]
\qquad , r \geq 1, 
$$
{\it symmetric} in the indices $i_1 \ldots i_r$. They are all
of the form $T \, (f \, U_{\psi}^*) = h \, f \, U_{\psi}^* $,
with $h \in  C^{\infty}€ (F \semi \G)$; in particular, they pairwise
commute.

A first observation is that the linear space
$$ {\frak h} (n) = \sum \Cb \cdot Y_i^{j} \oplus \sum \Cb \cdot X_k \oplus
   \sum  \Cb \cdot \d_{ab,i_1 \ldots i_r}^c 
$$
forms a Lie algebra and furthermore, if we let 
$$ \Hc (n) =  \Uc  ({\frak h} (n))
$$
denote the corresponding enveloping algebra, then
$\Hc (n)$ acts on $\Ac$ satisfying a Leibniz rule of the form
\begin{equation} \label{pr}
    h(ab) = \sum \, h_{(0)} \, (a) \, h_{(1)} \, (b) 
    \qquad \fl \, a,b \in {\Ac} .
\end{equation}

A second observation is that the product rules 
(\ref{YR}), (\ref{XR}) and (\ref{dijR}) 
can be converted into {\it coproduct} rules
\begin{equation}
\matrix{
\D \, Y_i^j = Y_i^j \ot 1 + 1 \ot Y_i^j ,  \hfill \cr
\cr
\D X_i = X_i \ot 1 + 1 \ot X_i + \sum_{k}€  \d_{ji}^k \ot Y_k^j ,\cr
\cr
\D \d_{ij}^k = \d_{ij}^k \ot 1 + 1 \ot \d_{ij}^k .  \hfill \cr
} 
\end{equation}
Together with the requirement 
$$
\D [Z_1 , Z_2] = [\D Z_1 , \D Z_2] \qquad \fl \, Z_1, Z_2 \in {\frak h}(n) ,
$$
which is satisfied on generators because of the flatness of the
connection, they
uniquely determine a multiplicative coproduct 
\begin{equation} \label{cp}
\D : {\Hc}(n) \ra {\Hc}(n) \otimes {\Hc}(n).
\end{equation}
One can check that this coproduct is coassociative and also that it is
compatible with the Leibniz rule (\ref{pr}), 
in the sense that
\begin{equation} \label{HA}
 \matrix{
 \D h = \sum \, h_{(0)} \ot h_{(1)} \qquad {\rm iff} \cr
 \cr
   \D (h_1 \, h_2) = \D h_1 \cdot \D h_2 \qquad \fl \, h_j \in {\Hc}(n) ,
    \qquad \fl \, a,b \in {\Ac} . \cr
    }
\end{equation}
Simple computations show that there is a unique antiautomorphism
$S$ of ${\Hc}(n)$ such that 
$$
 S(Y_i^{j}) = - Y_i^{j}, \qquad  S(\d_{ab}^c) = - \d_{ab}^c ,  \qquad
S(X_a )= -X_a +\, \d_{ab}^c Y_c^{b} .
$$
Moreover, $S$ serves as antipode for the bialgebra ${\Hc}(n)$; we
should note though that $S^{2}€ \neq I $.

To summarize, one has:

\begin{proposition}[\cite{CM2}] The enveloping algebra ${\Hc}(n)$ of
the Lie algebra generated by the canonical action of $\Rb^{n}€ \semi
{\frak gl} (n, \Rb)$ on $\Ac (n)$
has a unique coproduct which turns it into a Hopf algebra such that
its tautological action on $\Ac (n)$ is a Hopf action. 
\end{proposition}

\begin{remark}
{\rm The Hopf algebra ${\Hc}(n)$ acts canonically on the crossed 
product algebra $\Ac_{F(M)}€ = C_{c}^{\infty}€€ (F(M) \semi \G)$, 
for any flat affine manifold $M$ with a pseudogroup $\G$ of
orientation preserving diffeomorphisms. Using Morita 
equivalence, one can always reduce the case of a general
manifold $M$ to the flat
case. The obstruction one encounters in trying to transfer
the action of ${\Hc}(n)$, via
the Morita equivalence data,
from the flattened version to a non-flat $M$ is 
exactly the curvature of the manifold $M$.}  The analysis of this
obstruction in the general context of actions of Hopf algebras
on algebras should provide
the correct generalization of the notion of
Riemannian curvature in
the framework of noncommutative geometry.
\end{remark}

There is a more revealing definition (\cite{CM2})
of the Hopf algebra  ${\Hc}(n)$, 
in terms of a bicrossproduct construction (cf. e.g. \cite{K})
whose origin, in the case of finite groups, can be traced 
to the work of G. I. Kac. In our case,
it leads to the interpretation of ${\Hc}(n)$ 
as a bicrossproduct of two
Hopf algebras, $\Uc_{a}(n)€$ and  $\Sc_{u}(n)€$,
canonically associated to the decomposition of the
diffeomorphism group as a set-theoretic product
$$ {\rm Diff}(\Rb^{n}€) = G_{a}€(n) \cdot G_{u}€(n) , 
$$
where $G_{u}€(n)$ is  
the group of diffeomorphisms of the form  
$\psi (x) = \, x + \, o(x)$. 
$\Uc_{a}(n)€$ is just the universal enveloping 
$\Uc ({\frak g}_{a}€(n))$ of the group $ G_{a}€(n)$ of affine motions
of $\Rb^{n}€$, with its natural Hopf structure.  
$\Hc_{u}(n)€$ is the Hopf algebra of polynomial functions on the
pro-nilpotent group of formal diffeomorphisms associated to $ G_{u}€(n)$.

The modular character of the affine group
$ \d = \Trace : {\frak g}_{a}€(n) \ra \Rb$ extends to a character
$ \d \in {\Hc}(n)^{*}€$. It can be readily
verified that the corresponding twisted antipode satisfies
the involution condition $\wt{S}^2 = I$.
It follows that the pair $(\d, 1)$ fulfills (\ref{ic}) and hence
forms a modular pair in involution.

The preceding
bicrossproduct interpretation allows  
to relate the cyclic cohomology of ${\Hc}(n)$ with respect to
$(\delta, 1)$ to the Gelfand-Fuchs cohomology \cite{GF} of
the infinite-dimensional 
Lie algebra ${\frak a}_n $ of formal vector fields on $\Rb^n$.

\begin{theorem}[\cite{CM2}]  There is a canonical Van Est-type map of
    complexes which induces an isomorphism    
\begin{equation} \label{GF}    
  \sum_{i \equiv * (2)}€ H^{i}€({ \frak a}_n ) \simeq    
    H  C_{{\rm per} \, (\delta, 1)}^{*}€ \, (\Hc (n)). 
\end{equation}    
\end{theorem}
\smallskip

We now return to the spectral triple $(\Ac_{F}€, \Hc_{F}€, D)$
associated to $(M, \G)$ with $M$ flat. Then we have the canonical
Hopf action $ {\Hc}(n) \otimes \Ac_{F}€ \ra \Ac_{F}€$. In addition,
the crossed product $\Ac_{F}€$ inherits a canonical trace 
$ \tau_{F}€ : \Ac_{F}€ \ra \Cb$, dual to the volume form $vol_{F}€$,
\begin{equation} \label{tf}
\tau_{F}€ \, (f \, U_{\psi}^*) = 0 \quad \hbox{if} \quad \psi \ne 1 \, 
\quad \hbox{and} \quad
\ \tau_{F}€(f) =
\int_F f \, vol_{F}€ \, . 
\end{equation}
Using the $\G$-invariance of $vol_{F}€$, it is easy to check 
that the trace
$ \tau_{F}€$ is $\d$--invariant under the action of $\Hc (n)$,
i.e. the property (\ref{ip}) holds. We therefore obtain
a characteristic map
$$
\g_{F}€^{*}€ :  H  C_{{\rm per} \, (\delta, 1)}^{*}€ \, (\Hc (n)) 
   \ra H  C_{\rm per}^{*}€ (\Ac_{F}€) ,
$$
which together with (\ref{GF}) gives rise to a new characteristic
homomorphism:
\begin{equation} \label{charf}
      \chi_{F}€  :  H^* ({ \frak a}_n ) \ra 
      H  C_{\rm per}^{*}€ (\Ac_{F}€)  .
\end{equation}
Passing to $SO(n)$-invariants, one obtains an induced
characteristic map from the relative Lie algebra cohomology,
\begin{equation} \label{charp}
      \chi_{P}€  : \sum_{i \equiv * (2) }€ H^{i + p}€ ({ \frak a}_n , SO(n))
      \ra  H  C_{\rm per}^{*}€ (\Ac_{P}€) ,
\end{equation}
which instead of $ \tau_{F}€$  involves the analogous trace
$ \tau_{P}€$ of $\Ac_{P}$ and where $p = {\rm dim} P $.

Let us assume that the action of $\G$ on
$M$ has no degenerate fixed point. Recall the local formula (\ref{phi})
for $ch_{*}€(R)$.
Using the built-in 
affine invariance of a hypoelliptic affine operator
$R$, one can show that
any cochain on ${\Ac}_P$ of the form,
$$
\vp (a^0 , \ldots , a^n) = \, \int \!\!\!\!\!\! - \, a^0 [R,a^1]^{(k_1)}
\ldots
[R , a^n]^{(k_n)} \, |R|^{- (n + \vert 2k \vert)} ,  \qquad
\fl \, 
a^j
\in {\Ac}_P
$$
can be written as a finite linear combination
$$
\vp (a^0 , \ldots , a^n) = \sum_{\a} \, \tau_{P}€ (a^0 \, h_1^{\a} (a^1) 
\ldots h_n^{\a} (a^n)) , \quad {\rm with} \quad
h_{i}€^{\a}€ \in \Hc (n),
$$
and therefore belongs to the range of the characteristic map $\chi_{P}€ $.
The structure of the cohomology ring $H^* ({ \frak a}_n , SO(n))$ 
is well-known (\cite{Go}). It is computed by
the cohomology of the finite-dimensional complex
$$
 \{  E (h_1 , h_3 , \ldots , h_m) \ot P (c_1 , \ldots ,c_n) \ , \ d \}
$$
where $E (h_1 , h_3 , \ldots , h_m)$ is the exterior algebra in the
generators
$h_i$ of dimension $2i-1$, with $m$ the largest odd integer less than
$n$ and $i$
odd, while $P (c_1 , \ldots ,c_n)$ is
the polynomial algebra in the generators $c_i$ of degree $2i$
truncated by
the ideal of elements of weight $>2n$. The coboundary $d$ is defined by,
$$
dh_i = c_i \ , \ i \ \hbox{odd} \ , \ dc_i = 0  \ \hbox{for all} \ i \, .
$$
In particular,
the Pontryagin classes $p_i = c_{2i}$ are non-trivial for all
$2i \leq n$.

The final outcome of the preceding discussion is the following
index theorem
for transversely hypoelliptic operators on foliations:

\begin{theorem}[\cite{CM2}] 
   For any hypoelliptic affine operator $R$ on $P$, there exists
   a characteristic class 
   $\Lc (R) \in \sum_{i \equiv * (2)}€ H^{i + n}€ ({ \frak a}_n , SO(n))$
   such that
   $$  ch_{*}€(R)  = \chi_{P}€ (\Lc (R) ) \, \in \,
   H  C_{\rm per}^{*}€ (\Ac_{P}€) .
   $$
\end{theorem}
\medskip

\begin{remark}
{\rm We conclude with the remark that similar considerations, leading to 
analogous results, can be implemented
for more specialized cases of transverse geometries, such as 
complex analytic and symplectic (or Hamiltonian). 
For example, in the {\it symplectic case}, which is the less obvious
of the two,
the transverse data consists of a symplectic manifold
$(M^{2n}, \, \omega)$ together with a pseudogroup $\G_{sp}€$ of 
local symplectomorphisms. The corresponding frame bundle is the principal 
$Sp (n, \Rb)$--bundle $F_{sp}€$ of symplectic frames.  
Its quotient mod $K$, 
where $K =  Sp (n, \Rb) \cap O_{2n} \simeq U(n)$, is the bundle 
$P_{sp}€$ whose fiber at
$x \in M$ consists of the almost complex structures on $T_x M$ 
compatible with $\omega_x $. It carries an intrinsic 
{\it para--K\"ahlerian structure}, obtained as follows.
The typical fibre of $P_{sp}€$
can be
identified with the noncompact Hermitian symmetric space 
$Sp (n, \Rb)/U(n)$ and,
as such, it inherits a canonical K\"ahler structure. 
This gives rise to a natural K\"ahler structure on the vertical subbundle
$\Vc$ of the tangent bundle $TP_{sp}€$. On the other hand, the normal
bundle $\Nc = T P_{sp}€/ \Vc \simeq TM$ possesses a tautological
K\"ahler structure. Indeed,
a point in $P_{sp}€$ is by definition
an almost complex structure and thus, 
together
with $\omega$, determines a ``moving'' K\"ahler structure.
The entire
construction is functorial with respect to local symplectomorphisms. 
One
can therefore define the symplectic analogue $Q_{sp}€$
of the hypoelliptic signature operator as a graded direct sum
$$
Q_{sp}€ = (\bar{\partial}_V^* \, \bar{\partial}_V - 
\bar{\partial}_V \, \bar{\partial}_V^*) \op (\bar{\partial}_H +
\bar{\partial}_H^*) \, ,
$$
where the horizontal $\bar{\partial}$-operator  
$\bar{\partial}_H$ is associated to a symplectic connection.}
The corresponding index theorem asserts that, 
with $(M^{2n}, \, \omega)$ flat,  acted upon by
an arbitrary pseudogroup of local symplectomorphisms $\G_{sp}$,
and with
\begin{equation} \label{charsp}
      \chi_{sp}€  :  H^* ({\frak a}^{sp}€_n , \, U(n)) \ra 
      H  C_{\rm per}^{*}€ (\Ac_{P_{sp}€}€) ,
\end{equation}
denoting the 
characteristic map corresponding to
the Lie algebra ${\frak a}^{sp}€_n$
of formal Hamiltonian vector fields on $\Rb^{2n}€$,
 there exists
a characteristic class $\Lc_{sp}€ \in  
\sum_{i \equiv * (2)}€ H^{i}€ ({\frak a}^{sp}€_n , \, U(n))$
such that 
$$  ch_{*}€(Q_{sp}€)  = \chi_{sp}€ (\Lc_{sp}€ ) \, \in \,
   H  C_{\rm per}^{*}€ (\Ac_{P_{sp}}€) .
$$
{\rm As in the para--Riemannian case, the proof relies on the local Chern
character formula (\ref{phi}) and on the cyclic cohomology of
the Hopf algebra $\Hc_{sp}€(n)$, associated to the group of
symplectomorphism of $\Rb^{2n}€$ in the same manner as $\Hc (n)$ was
constructed from ${\rm Diff} (\Rb^{n}€)$.}
\end{remark}
\vspace{1cm}

\section{Quantum groups and the modular square}

The definition of cyclic co/homology of Hopf algebras hinges on the
existence of modular pairs in involution. The necessity of 
this condition may appear as artificial. 
In fact, quite the opposite is true and the
examples given below serve to illustrate that most Hopf algebras
arising in ``nature'', including
quantum groups and their duals,
do come equipped with intrinsic modular
pairs.
\medskip

{\bf 1.} We begin with the class of {\it quasitriangular Hopf
algebras}, introduced by Drinfeld \cite{Dr},
in connection with the quantum inverse scattering method for
constructing quantum integrable systems. 
Such a Hopf algebra comes endowed 
with an universal $\Rc$--{\it matrix}, inducing solutions of the
Yang-Baxter equation on each of their modules. (For
a lucid introduction into the subject, see \cite{K}).

A {\it quasitriangular Hopf algebra} is  a Hopf algebra
$\Hc$ which admits 
an invertible element 
$R = \sum_{i}€ s_{i}€ \ot t_{i}€  \in \Hc \ot \Hc $,
such that
\begin{equation} \label{Rm}
 \matrix{
 \D^{\rm op}€(x) = R \D(x) R^{-1}€ , \quad  \fl x \in \Hc \, , \cr
 \cr
  (\D \ot I )(R) = R_{13}€R_{23}€  \, , \cr
  \cr
  (I \ot \D) (R) =  R_{13}€R_{12}€   ,\,  \cr
    }
\end{equation}
where we have used the customary ``leg numbering'' notation, e.g.
$$ R_{23}€ = \sum_{i}€ 1 \ot s_{i}€ \ot t_{i}€ .
$$
The square of the antipode is then an 
inner automorphism,
$$
S^2 (x) = u x u^{-1} \, ,
$$
with
$$
u = \sum \, S ( t_{i}€) s_{i}€ , quad
u^{-1}€ = \sum \, S^{-1}€ (t_{i}€)  S(s_){i}€ .
$$
Furthermore, $ u S(u) = S(u) u $ is central in $\Hc$ and one has
$$ \varepsilon (u) = 1, \quad
\D u = (R_{21} R)^{-1} (u \otimes u) =  (u \otimes u) (R_{21} R)^{-1} \, .
$$

A quasitriangular Hopf algebra $\Hc$ is called a {\it ribbon algebra} 
\cite{RT}, if there exists a {\it central element} $\theta \in \Hc$ 
such that
\begin{equation} \label{rib}
\Delta (\theta) = (R_{21} \, R)^{-1} \, (\theta \otimes \theta) \ , \quad 
\varepsilon (\theta) = 1 \ , \quad S(\theta) = \theta \, .
\end{equation}
Any  quasitriangular Hopf $\Hc$ algebra has 
a ``double cover'' cover ${\wt \Hc}$ satisfying the ribbon condition 
(\ref{rib}). More precisely (cf.\cite{RT}), 
$$
\wt \Hc = \Hc \, [\theta] / (\theta^2 - u \, S(u))
$$
has a unique Hopf algebra structure such that, under the natural 
inclusion, $\Hc$ is a Hopf subalgebra. 

If $\Hc$ is a ribbon algebra, by setting
$$
\sigma = \theta^{-1} \, u \, ,
$$
one gets a group-like element
$$
\Delta \, \sigma = \sigma \otimes \sigma \ , \quad \varepsilon (\sigma) = 1 \ , 
\quad S(\sigma) = \sigma^{-1} \, 
$$
such that, for any $x \in \wt \Hc$,
$$
\matrix{
({\s^{-1}€\circ S)}^2 (x) &= &\sigma^{-1} S (\sigma^{-1}  S (x)) = 
\sigma^{-1} \, S^2 (x) \sigma  \hfill \cr
&= &\sigma^{-1} u x u^{-1} \sigma = \theta  x  \theta^{-1} = x 
\, . \hfill \cr
}
$$
Thus, $(\ve, \sigma)$ is a modular pair in involution
for $\Hc$.
\smallskip

By dualizing the above definitions one obtains the notion of 
{\it coquasitriangular}, resp.  {\it coribbon algebra}. Among the most
prominent examples of coribbon algebras are the function
algebras of the classical quantum groups  ($GL_q (N)$, $SL_q (N)$,
$SO_q (N)$, $O_q (N)$ and $Sp_q (N)$). 
The analogue of
the above {\it ribbon group-like element} $\sigma$
for a coribbon algebra $\Hc$, is
the {\it ribbon character} $\delta \in \Hc^*$. The corresponding 
twisted antipode 
satisfies the condition ${\wt S}^2 = 1$, which renders $(\delta, 1)$
as a canonical modular pair in involution for $\Hc$.

We thus have: 
\smallskip

\begin{proposition}[\cite{CM3}] Coribbon algebras and compact quantum
    groups are each intrinsically  endowed with a modular pair in 
    involution $(\d, 1)$. Dually, ribbon algebras and duals of
    compact quantum groups are each intrinsically endowed with a
    modular pair in involution $(1, \s)$. 
\end{proposition}

For a {\it compact quantum
group} in the sense of Woronowicz, the stated property follows from
Theorem 5.6 of \cite{W}, describing the
modular properties of the analogue of Haar measure.
\medskip

{\bf 2.} Evidently, one can produce modular pairs in involution
with both $\d$ and $\s$ nontrivial by forming tensor products of
dual classes of Hopf algebras as in the preceding statement.  
The fully non-unimodular situation arises
naturally however, in the case of
{\it locally compact quantum groups}, because of the existence,
by fiat or otherwise, of left and right Haar weights.
We refer to \cite{KV} for the most recent and concise 
formalization of this notion, which is in remarkable agreement
with
our framework for cyclic co/homology of Hopf algebras
and of Hopf actions. This accord is manifest in the 
following {\it construction
of a  modular square} associated to a Hopf algebra $\Hc$
modelling a locally compact group.
Since the inherent analytic intricacies
are beyond the scope of
the present exposition, we shall keep the illustration at a
formal level (comp. \cite{V} for an algebraic setting). 
\smallskip

By analogy with the definition of a $C^{*}€$--algebraic quantum
group in \cite{KV}, we assume the existence and uniqueness (up to
a scalar)
of a {\it left invariant
weight} $\vp$, satisfying a
{\it KMS--like condition}. The invariance means that
\begin{equation} \label{linv}
 (I \ot \vp) ((I \ot x) \D (y)) = S((I \ot \vp) \D (x) (I \ot y))  ,
 \qquad \fl \, x, y \in \Hc ,
\end{equation}
while the KMS condition stipulates the existence of a 
{\it modular group
of automorphisms} $\s_{t}€$ of $\Hc$, such that
\begin{equation} \label{lkms}
   \vp \circ \s_{t}€ = \vp \,  , \qquad  \vp (xy) = 
   \vp (\s_{i}€ (y) x) , \quad \fl \, x, y \in \Hc \, .
\end{equation}
Taking $\psi = \vp \circ S^{-1}€$, one obtains a 
a {\it right invariant weight}, 
\begin{equation} \label{rinv}
  (\psi \ot I) (\D(x) (y \ot 1)) = S( (\psi \ot I) (x \ot 1) 
  \D (y))  ,
\end{equation}
which is also unique up to a scalar and has
modular group $\s^{\prime}€_{t}€ = S \circ \s_{t}€ \circ S^{-1}€$.

It will be convenient to express the above properties in terms of
the natural left and right actions of the dual Hopf algebra
$\hat \Hc = \Hc^{*}€$ on $\Hc$.
Given $\om \in \hat \Hc$ and $x \in \Hc$,
we denote
\begin{equation} \label{act}
 \matrix{
   \om \cdot x = \sum \om (x_{(1)}€) \, x_{(2)}€ = 
   (\om \ot I) \D(x) \, ,\cr
 \cr
  x \cdot \om = \sum   x_{(1)}€ \, \om (x_{(2)}€)  = 
   (I \ot \om) \D(x) \, . \cr
    }
\end{equation}
With respect to the natural product of $\hat \Hc$,
$$ ( \om_{1}€ \ast \om_{2}€) (x)  =  <  \om_{1}€ \ot \om_{2}€ , \, \D(x) >
  =  \sum  \om_{1} (x_{(1)}) \, \om_{2}€ (x_{(2)}) \, , \qquad
  \fl x \in \Hc,
$$
the left action in (\ref{act}) is the transpose of the left regular
representation of $\hat \Hc$, hence defines a representation
of the opposite algebra ${\hat \Hc}^{\rm op}€$ on $\Hc$, while 
the right action in (\ref{act}), being the transpose of the
right regular representation of $\hat \Hc$, gives a representation
of the algebra $\hat \Hc$ on $\Hc$.  On the other hand, it is
easy to check that both 
actions satisfy the rule (\ref{HL}) and therefore (\ref{act})
{\it defines a Hopf action of the tensor product Hopf algebra 
$ \wt \Hc := {\hat \Hc}^{\rm op} \ot \hat \Hc$ 
on $\Hc$}},
\begin{equation} \label{hac}
   \wt \Hc \ot  \Hc \ra \Hc \, , \quad
   (\om_{1}€ \ot \om_{2}€ , \, x) \ra \om_{1}€ \cdot x \cdot \om_{2}€  \, .
\end{equation}   
The invariance conditions (\ref{linv}) and (\ref{rinv}) 
can now be rewritten as
follows:
\begin{equation} \label{inv}
 \matrix{
  \vp (( \om \cdot x) \, y) = \vp (x \, (S^{-1}€(\om) \cdot y))  , \cr
 \cr
 \psi ((x \cdot \om) \, y) = \psi (x \, (y \cdot S(\om))) , \cr
    }
\end{equation}
where $S^{-1}$ occurs in the first identity as the antipode of
${\hat \Hc}^{\rm op}$.

The left invariance property of $\psi$ gives the analogue
of the modular function, namely a {\it  group-like element}
$\d \in \Hc$ such that
$$ (I \ot \psi) \D (x) = \psi (x) \d \, , \qquad \fl \, x \in \Hc,
$$
or equivalently
\begin{equation} \label{del}
    \psi (\om \cdot x) = \om (\d) \psi (x) \, , \qquad \fl \, x \in 
    \Hc.
\end{equation}
The {\it modular element} $\d$ also relates
the left and right Haar weights: 
\begin{equation}
\vp (x) = \psi (\d^{1 \over 2}€ \, x \, \d^{1 \over 2}€) 
  \qquad \fl \, x \in \Hc.
\end{equation}
In particular, the full invariance property of $\vp$ under the action of
$ \wt \Hc$ is given by
\begin{equation} \label{fli}
    \matrix{
    \vp ((\om_{1}€ \cdot x \cdot \om_{2}€) \, y) =
    \vp (x \,  (S^{-1}€(\om_{1}€) \cdot y \cdot S_{\d^{-1}€}€ 
    (\om_{2}€))) \, ,\cr
    \cr
    \fl \quad  \om_{1}€ , \om_{2}€  \in \wt \Hc , \quad x, y \in \Hc \, ,
    }
\end{equation}    
where $ S_{\d^{-1}}$ denotes the twisted antipode (\ref{dS})
corresponding to $\d^{-1}€$.
\smallskip

Let us now form the {\it midweight} $\tau$,
\begin{equation} \label{mid}
    \tau (x) = \vp ( \d^{ - {1 \over 4}}€ \, x \,  \d^{ - {1 \over 4}}€)
    \qquad \fl \, x \in \Hc .
\end{equation} 
One checks that its behavior under the action of $ \wt \Hc$ is as 
follows:
\begin{equation} \label{fmi}
    \tau ((\om_{1}€ \cdot x \cdot \om_{2}€) \, y) =
    \tau (x \,  (S^{-1}€_{\d^{1/2}€}€
    (\om_{1}€) \cdot y \cdot S_{\d^{-1/2}€}€ 
    (\om_{2}€))) \, ,
\end{equation}    
for any $ \om_{1}€ , \om_{2}€  \in \wt \Hc $ and $x, y \in \Hc$.
In other words, $\tau$ is $\wt \d$--{\it invariant under the action 
(\ref{hac}), with respect to the character}
\begin{equation} \label{wdel} 
\wt \d = \d^{1 \over 2}€ \ot \d^{- {1 \over 2}}€ \in 
      {\wt \Hc}^{\ast}€ = \Hc^{\rm op}€ \ot \Hc .
\end{equation}
On the other hand, it follows as in \cite{KV}, 
but with a slightly different notation, that 
there exists a group-like element $\s \in \hat \Hc$, such that
the modular groups of $\vp , \psi$ can be expressed as follows:
\begin{equation} \label{lrg}
    \matrix{
   \s_{t}€(x) = 
   \d^{it/2}€(\s^{it/2}€ \cdot x \cdot \s^{it/2}€) 
   \d^{-it/2}€ \, , \cr
\cr
 \s^{\prime}€_{t}€(x) = 
   \d^{-it/2}€(\s^{it/2}€ \cdot x \cdot \s^{it/2}€) 
   \d^{it/2}€ \, , \quad \fl \, x\in \Hc .
}
\end{equation}
In terms of the modular group of $\tau$, to be denoted 
$\s^{\tau}_{t}$ €,
(\ref{lrg}) is equivalent to 
\begin{equation} \label{mg}
   \s^{\tau}_{t} (x) = \s^{it/2}€ \cdot x \cdot \s^{it/2}€ \, ,
    \quad \fl \, x\in \Hc .
\end{equation} 
This shows that $\tau$ is a $\wt \s$--{\it trace for the 
action (\ref{hac})  ,
with group-like element}
\begin{equation} \label{wg} 
\wt \s = \s^{1 \over 2}€ \ot \s^{1 \over 2}€ \in 
      {\wt \Hc}^{\ast}€ = \Hc^{\rm op}€ \ot \Hc .
\end{equation}
It remains to compute the square of the corresponding 
doubly twisted antipode of $\wt \Hc$,
\begin{equation} \label{dt}
{\wt \s}^{-1}€ \circ  S_{\wt d}€ = 
\s^{-1/2}€ \circ S^{-1}€_{\d^{1/2}€}€ \ot 
\s^{-1/2}€ \circ S_{\d^{-1/2}€}€ : {\hat \Hc}^{\rm op} \ot \hat \Hc
  \ra {\hat \Hc}^{\rm op} \ot \hat \Hc\, .
\end{equation}
It suffices to compute the square of
$ \s^{-1/2}€ \circ S_{\d^{-1/2}€}€ : {\hat \Hc} \ra {\hat \Hc} \, ,$
or equivalently, the square of its transpose, for which a straightforward
calculation gives:
$$ ( \s^{-1/2}€ \circ S_{\d^{-1/2}€}€)^{2}€ = 
   <\s^{-1/2}€, \d^{-1/2}€ >\, I_{\hat \Hc}€ \, .
$$
Since the passage to the opposite algebra gives the reciprocal scalar,
it follows that
\begin{equation} \label{sq}
( {\wt \s}^{-1}€ \circ  S_{\wt d}€ )^{2}€ = I_{\wt \Hc}€ \, .
\end{equation}
\smallskip

We summarize the conclusions of the preceding discussion  
in the following result:
\smallskip

\begin{theorem} 
    {\rm (i)} The Hopf algebra $\wt \Hc = {\hat \Hc}^{\rm op} \ot \hat \Hc$ 
         possesses a canonical modular pair in involution 
         $ (\wt \d = \d^{1 \over 2}€ \ot \d^{- {1 \over 2}}€  , \,
	 \wt \s = \s^{1 \over 2}€ \ot \s^{1 \over 2}€)$. 
 
    {\rm (ii)} The Haar midweight $\tau$, given by (\ref{mid}),
        is a $\wt \d$--invariant $\wt \s$--trace for the canonical
        action of $\wt \Hc$ on $\Hc$. 
\end{theorem}
\smallskip

\noindent The first statement characterizes the construction
we referred to as {\it the modular square} associated to
a Hopf algebra $\Hc$ that models a locally compact quantum group.
Together with the Haar midweight $\tau$ of the second
statement, it determines
in cyclic cohomology a {\it modular characteristic homomorphism} 
\begin{equation} \label{mch}
     \g_{\tau}€^{*}€ : H  C_{({\wt \d}, {\wt \s})}^* ({\wt \Hc}) \ra
     H  C^* (\Hc) \, .
\end{equation}


\vspace{1cm}
 


\begin{thebibliography}{99}
        
\bibitem{Ad} Adams, J. F., On the cobar construction, 
\textit{Proc. Nat. Acad. Sci. USA} \textbf{42} (1956), 409-412.

\bibitem{BG} Beals, R. and Greiner, P., 
\textbf{Calculus on Heisenberg manifolds}, 
Annals of Math. Studies {\bf 119},
Princeton Univ. Press, Princeton, N.J., 1988.


\bibitem{Bu} D. Burghelea, The cyclic homology of the group rings, 
\textit{Comment. Math. Helv.} {\bf 60} (1985), 354-365.


\bibitem{Ca} Cartier, P.,  Cohomologie des coalg\`ebres (Expos\'e 5). 
In \textit{S\'eminaire Sophus Lie, 1955-56}, Facult\'e des Sciences de 
Paris, 1957.

\bibitem{C0} Connes, A., Spectral sequence and homology of currents for 
operator algebras, \textit{Math. Forschungsinst. Oberwolfach Tagungsber.} 
41/81; \textit{ Funktionalanalysis und $C^*$-Algebren}, 27-9/3-10, 1981.
    
\bibitem{C1} Connes, A., $C^*$ alg\`ebres et g\'eom\'etrie
differentielle, \textit{C.R. Acad. Sci. Paris}, Ser.~A-B \textbf{290} 
(1980), 599-604.

\bibitem{Co0} A. Connes, Noncommutative differential geometry. Part I: 
The Chern character in $K$-homology, 
\textit {Preprint IHES} (M/82/53), 1982; Part 
II: de Rham homology and noncommutative algebra, 
\textit { Preprint IHES} (M/83/19), 1983.
       
\bibitem{C2} Connes, A., Noncommutative differential geometry, 
\textit{Inst. Hautes Etudes Sci. Publ. Math.} \textbf{62} (1985), 
257-360.

\bibitem{C3} Connes, A.,  Cohomologie cyclique et foncteur $Ext^n$, 
\textit{C.R. Acad. Sci. Paris}, Ser.I Math \textbf{296} (1983),
953-958.
   
\bibitem{C4} Connes, A., Cyclic cohomology and the
transverse fundamental class of a foliation. In \textit{ Geometric
methods in operator algebras}, (Kyoto, 1983), pp. 52-144,
\textit{ Pitman Res. Notes in Math.} \textbf{123} Longman,
Harlow, 1986.

\bibitem{Co} Connes, A., \textbf{Noncommutative geometry}, Academic
Press, 1994.

\bibitem{CM0} Connes, A. and Moscovici, H., Cyclic
cohomology, the Novikov conjecture and hyperbolic groups,
\textit{Topology} \textbf{29} (1990), 345-388.

\bibitem{CM1}  Connes, A. and Moscovici, H., The local index
formula in noncommutative geometry, \textit{GAFA}  \textbf{5} (1995),
174-243.

\bibitem{CM2} Connes, A. and Moscovici, H., Hopf Algebras, Cyclic 
Cohomology and the Transverse Index Theorem, \textit{Commun. Math. 
Phys.} \textbf{198} (1998), 199-246 . 

\bibitem{CM3} Connes, A. and Moscovici, H.,  Cyclic 
Cohomology and Hopf Algebras, \textit{Letters Math. Phys.} 
\textbf{48}  (1999), 97-108.

\bibitem{Cr} Crainic, M., Cyclic Cohomology of Hopf Algebras and a 
Noncommutative Chern-Weil Theory, \textit{Preprint} QA/9812113.


\bibitem{Dr} Drinfeld, V. G., Quantum groups. In \textit{Proc.
Int. Cong. Math.} (Berkeley, 1986), pp. 798-820, Amer. Math. Soc.,
Providence, RI, 1987.


\bibitem{GF} Gelfand, I. M. and Fuchs, D. B., Cohomology of the 
Lie algebra of formal vector fields, 
\textit{Izv. Akad. Nauk SSSR} \textbf{34} (1970), 322-337.

\bibitem{Go} Godbillon, G., Cohomologies d'alg\`ebres de Lie de 
champs de vecteurs formels, 
\textit{Seminaire Bourbaki (1971/1972), Expos\'e No.421,}
Lecture Notes in Math. Vol. \textbf{383},
Springer, Berlin (1974), 69-87.

\bibitem{HS} Hilsum, M. and Skandalis, G., : Morphismes K-orient\'es 
d'espaces de feuilles et fonctorialit\'e en th\'eoriede Kasparov, 
\textit {Ann. Sci. Ecole Norm. Sup. (4)} \textbf{20} (1987), 325-390.

\bibitem{K} Kassel, C., \textbf{Quantum Groups}, Springer-Verlag, 
1995.

\bibitem{KV} Kustermans, J. and Vaes, S., A simple definition for 
locally compact quantum groups, 
\textit{C.R. Acad. Sci. Paris}, Ser.I Math \textbf{328} (1999),
871-876.

\bibitem{L} Loday, J.L., \textbf{Cyclic Homology}, Springer-Verlag, 
1992, 1998.

\bibitem{RT} Reshetikhin, N. Yu. and Turaev, V. G., Ribbon graphs and 
their invariants derived from quantum groups,
\textit{Commun. Math. Phys.} \textbf{127} (1990), 1-26 .

\bibitem{S} Sweedler, M.E., \textbf{Hopf Algebras}, W.A.~Benjamin,
Inc., New York, 1969.

\bibitem{Ts} Tsygan, B. L., Homology of matrix Lie algebras over rings and 
the Hochschild homology, \textit{ Uspekhi Math. Nauk.}
\textbf{38} (1983), 217-218.

\bibitem{V} Van Daele, A., An algebraic framework for group duality, 
\textit{Adv. Math.} \textbf{140} (1998), 323-366.

\bibitem{W} Woronowicz, S.L., Compact matrix pseudogroups, 
\textit{Commun. Math. Phys.} \textbf{111} (1987), 613-665.

\end{thebibliography}
\end{document}